\documentclass{icmart}

\contact[speicher@math.uni-sb.de]{Saarland University, Fachrichtung Mathematik, Postfach151150, 66041 Saarbr\"ucken, Germany}

\newtheorem{theorem}{Theorem}[section]

\theoremstyle{definition}
\newtheorem{definition}[theorem]{Definition}
\newtheorem{remark}[theorem]{Remark}

\newtheorem{example}[theorem]{Example}

\def\la{\langle}
\def\ra{\rangle}
\def\id{\text{{id}}}

\def\cA{\mathcal{A}}

\def\cP{\mathcal{P}}
\def\cB{\mathcal{B}}

\def\cX{\mathcal{X}}

\def\CC{{\Bbb C}}
\def\NN{{\Bbb N}}
\def\RR{{\Bbb R}}

\def\ff{\varphi}
\def\tr{\text{tr}}
\def\ee{\varepsilon}

\def\kk{\kappa}

\title[Free probability and random matrices]{Free probability and random matrices}

\author[Roland Speicher]
{Roland Speicher\thanks{
The author was partially supported by ERC Advanced Grant NCDFP 339760 during the preparation of this article.}}
\begin{document}

\begin{abstract}
The concept of freeness was introduced by Voiculescu in the context of
operator algebras. Later it was observed that it is also relevant for large random matrices. We will show how the combination of various free probability results with a linearization trick allows to address successfully the problem of determining the asymptotic eigenvalue distribution of general selfadjoint polynomials in independent random matrices. 
\end{abstract}



\maketitle

\section{Introduction}
Free probability theory was introduced by Voiculescu around 1983 in order to attack the isomorphism problem of von Neumann algebras of free groups. Voiculescu isolated a structure showing up in this context which he named ``freeness''.
His fundamental insight was to separate this concept from its operator algebraic origin and investigate it for its own sake. Furthermore, he promoted the point of view that freeness should be seen as an (though non-commutative) analogue of the classical probabilistic concept of ``independence'' for random variables. Hence freeness is also called ``free independence'' and the whole subject became to be known as ``free probability theory''.

The theory was lifted to a new level when Voiculescu discovered in 1991 that the freeness property is also present for many classes of random matrices, in the asymptotic regime when the size of the matrices tends to infinity. This insight, bringing together the apriori entirely different theories of operator algebras and of random matrices, had quite some impact in both directions. Modelling operator algebras by random matrices resulted in some of the deepest results about operator algebras of the last decades; whereas tools developed in operator algebras and free probability theory could now be applied to random matrix problems, yielding in particular new ways to calculate the asymptotic eigenvalue distribution of many random matrices. Since random matrices are also widely used in applied fields, like wireless communications or statistics, free probability is now also quite common in those subjects. 

Whereas Voiculescu's original approach to free probability is quite analytic and operator algebraic in nature, I provided another, more combinatorial, approach. This rests on the notion of ``free cumulants'' and is intimately connected with the lattice of ``non-crossing partitions''.

In this lecture we will give an introduction to free probability theory, with focus on its random matrix and combinatorial side. Freeness will be motivated not by its initial occurrence in operator algebras, but by its random matrix connection.
The main result we are aiming at is also a very general random matrix problem, namely how to calculate the distribution of selfadjoint polynomials in independent random matrices. Whereas there exists a vast amount of literature on solving this problem for various special cases, often in an ad hoc way, we will see that free probability gives a conceptual way to attack this problem in full generality. 

For more information on other aspects of the subject one might consult the earlier ICM contributions of Voiculescu \cite{V-ICM}, Haagerup \cite{H-ICM}, and Shlyakhtenko \cite{S-ICM} (for the operator algebraic aspects of free probability) or of Biane \cite{B-ICM} (for applications to the asymptotics of representations of symmetric groups). Extensions of the theory to rectangular matrices can be found in \cite{BGe}, and to ``second order freeness'' (describing fluctuations of random matrices) in \cite{CMSS}.
The monographs \cite{HP, MSp, NS, VDN} give general introductions to free probability; \cite{NS} has its main emphasis on the combinatorial side of the subject.
The applied side of free probability is addressed, for example, in \cite{CB,RE,TV}. 

\section{Motivation of freeness via random matrices}
 
In this chapter we want to motivate the definition of freeness on the basis of random matrices.

\subsection{Asymptotic eigenvalue distribution of random matrices} 
We are interested in computing the eigenvalue distribution of $N\times N$ random matrices as $N\rightarrow\infty$. Here and in the following we will only consider selfadjoint random matrices. This guarantees that the eigenvalues are real and strong analytical tools are available to deal with such situations. For non-selfadjoint matrices the eigenvalues are in general complex and the situation, especially in the case of non-normal matrices, is more complicated. We will make some remarks on this situation at the very end of our lecture, in Sect.~\ref{sect:outlook}.

The typical feature for many basic random matrix ensembles is the almost sure convergence to some limiting eigenvalue distribution. Furthermore, quite  often this limiting distribution can be effectively calculated.

\begin{example}
We consider an $N\times N$ \emph{Gaussian random matrix}. This is a selfadjoint matrix $X_N=\frac{1}{\sqrt{N}}(x_{ij})_{i,j=1}^N$
such that the entries $\{ x_{ij}\} _{i\geq j}$ are independent and identically distributed complex (real for $i=j$) Gaussian random variables with mean 
${E}[x_{ij}]=0$ and variance ${E}[x_{ij}\bar x_{ij}]=1$.

The following figure shows typical histograms of the eigenvalues of Gaussian random matrices, for different values of $N$.

\begin{figure}[h]
\includegraphics[width=1.6in]{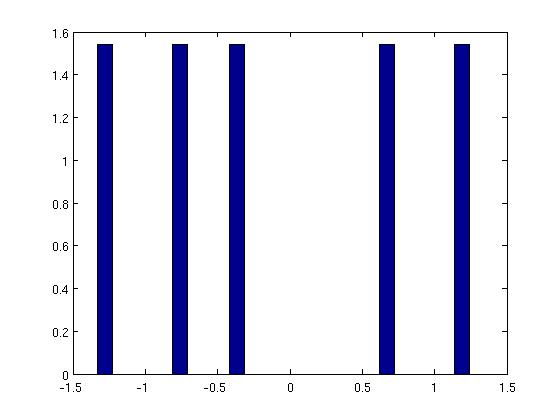}
\includegraphics[width=1.6in]{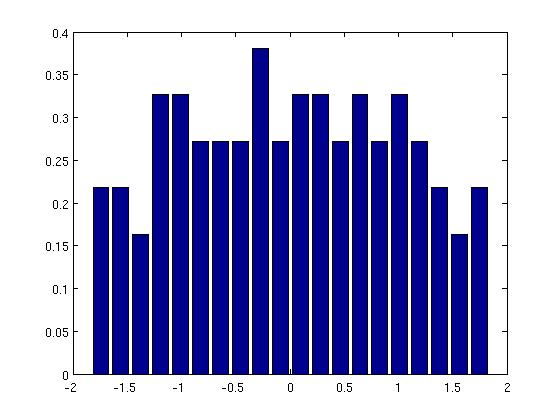}
\includegraphics[width=1.6in]{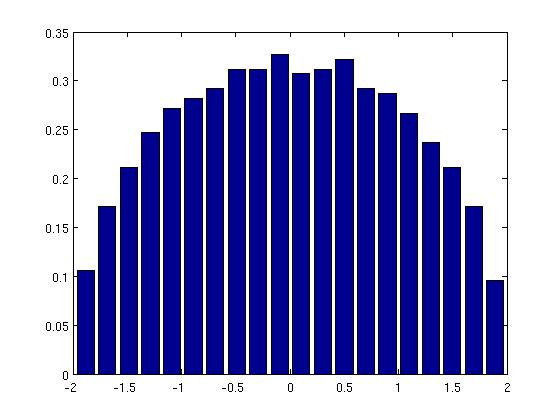}
\caption{\label{fig:eins} Histogram of the $N$ eigenvalues for a realization of an $N\times N$ Gaussian random matrix; for N= 5, 100, 1000}
\end{figure}

One sees that, whereas for small $N$ there is no clear structure, for large $N$ the eigenvalue histogram is approaching a smooth curve. Actually, this curve is deterministic, it is (almost surely) always the same, independent of the actual realized matrix from the ensemble. What we see here, is one of the first and most famous
results of random matrix theory: 
for such matrices we have almost sure convergence to \emph{Wigner's semicircle law}, given by the densitiy $\rho (t)=\frac{1}{2\pi}\sqrt{4-t^2}$.
In Fig.~\ref{fig:semi} we compare one realization of a $4000\times 4000$
Gaussian random matrix with the semicircle.
\end{example}  

\begin{example}
An other important class of random matrices are\emph{Wishart matrices}; those are of the form $X_N=A_NA_N^*$, where $A_N$ is an $N\times M$ matrix with independent Gaussian entries. If we keep the ratio $M/N$ fixed,
its eigenvalue distribution converges almost surely to the \emph{Marchenko-Pastur distribution}; see Fig.~\ref{fig:semi}.

\begin{figure}[h]
\includegraphics[width=2.6in]{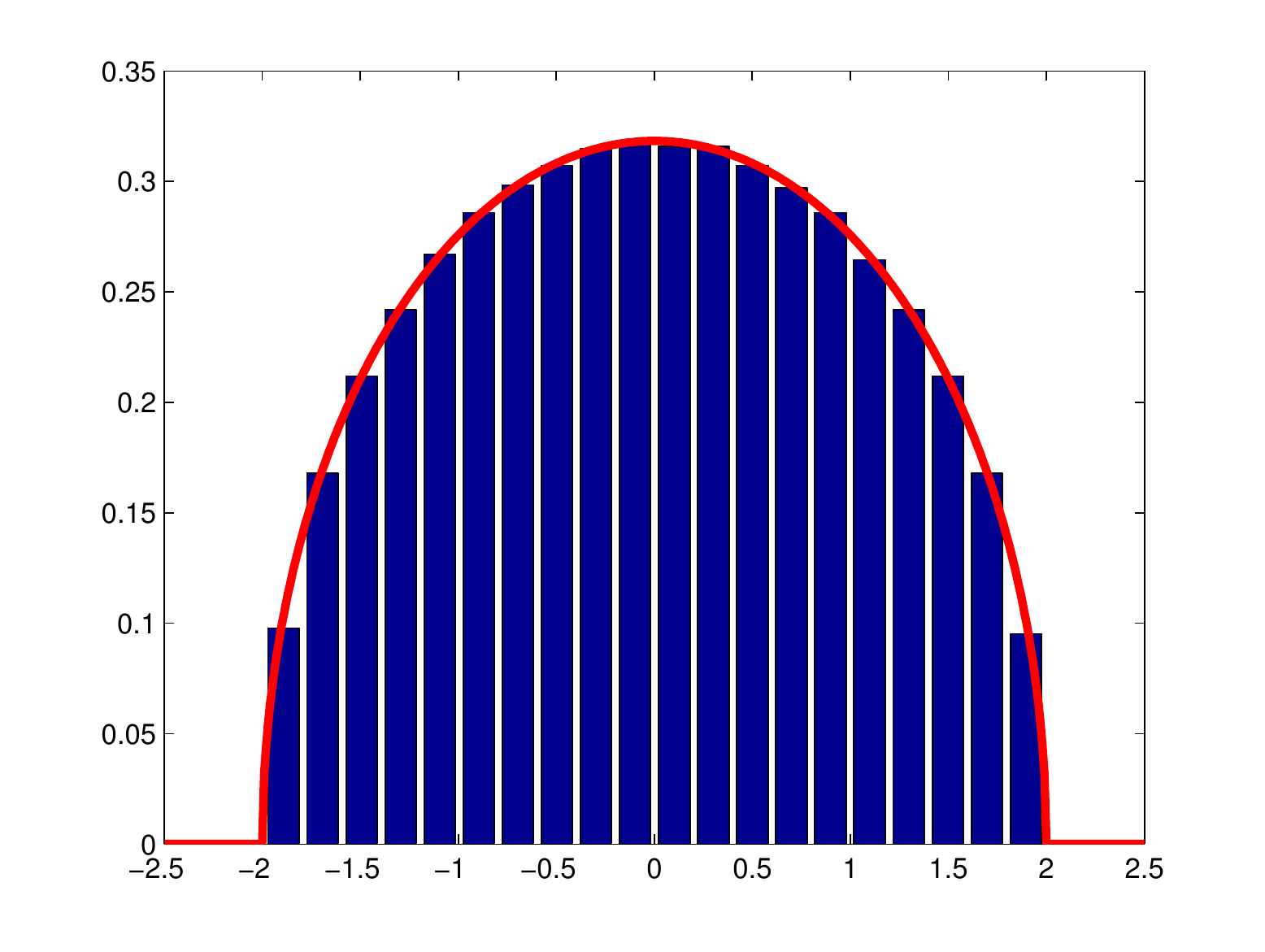}
\includegraphics[width=2.6in]{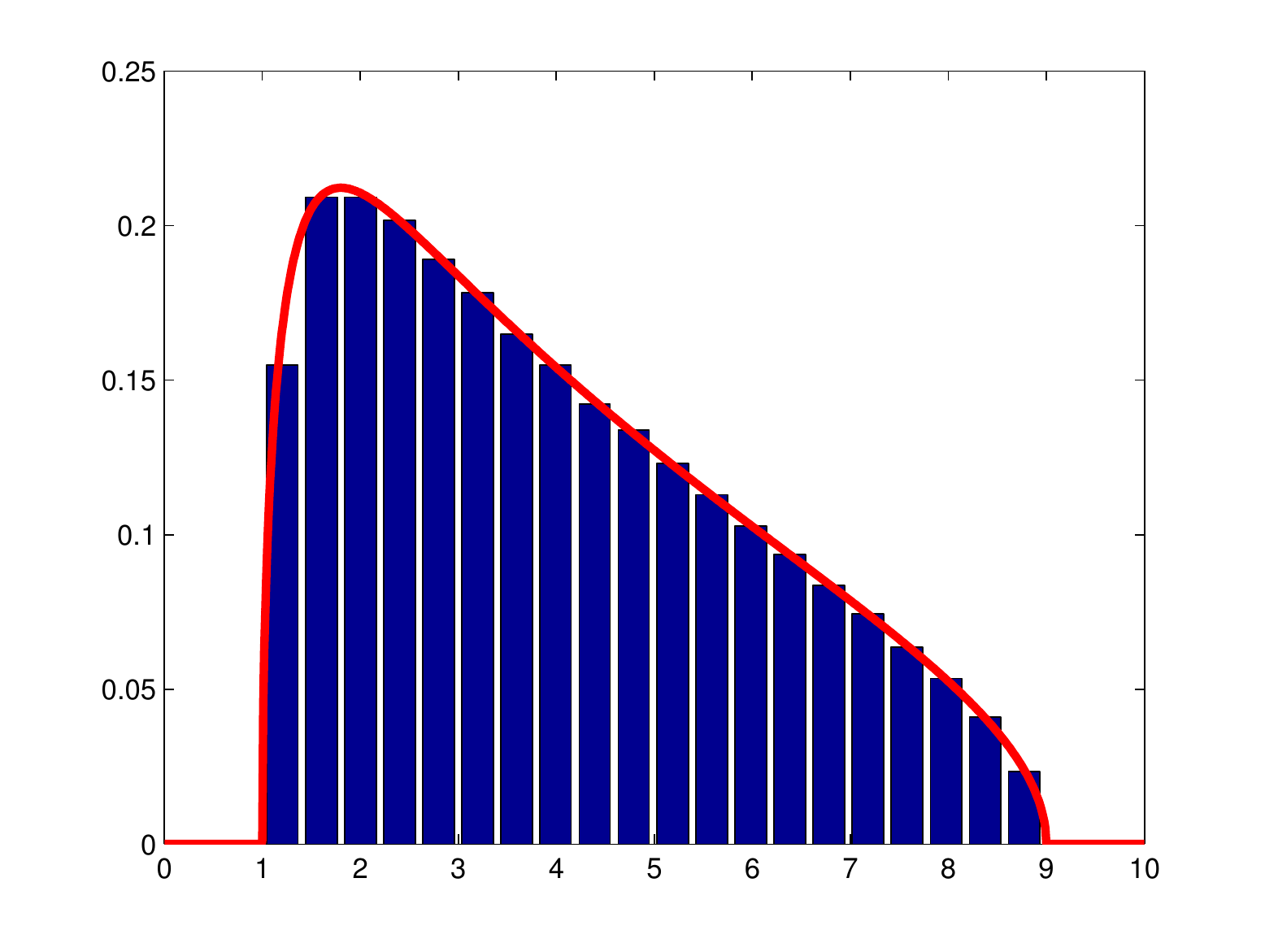}
\caption{\label{fig:semi} Comparision between the histogram of one realization of a random matrix and the analytic form of the density in the $N\to\infty$ limit; left: $4000\times 4000$ Gaussian random matrix versus semicircle distribution; right: Wishart random matrix with $N=3000$, $M=12000$ versus the corresponding Marchenko-Pastur distribution}
\end{figure}
\end{example}

\subsection{Polynomials in several random matrices}

Instead of looking at one-matrix ensembles we are now interested in the case of several matrices.  Let us consider two sequences $X_N$ and $Y_N$ of
selfadjoint $N\times N$ matrices such that both sequences have an asymptotic eigenvalue
distribution for $N\to\infty$. We are interested in the asymptotic eigenvalue
distribution of sequences $p(X_N,Y_N)$ for some non-trivial functions $p$
of two non-commuting variables.  We will restrict to the simplest class of functions, namely $p$ will be a (non-commutative) polynomial. 
As mentioned before, we are only dealing with selfadjoint matrices, thus $p$ should be a selfadjoint polynomial in order to ensure that also $p(X_N,Y_N)$ is selfadjoint.

In general, the distribution of $p(X_N,Y_N)$ will depend on the relation between the eigenspaces of $X_N$ and of
$Y_N$. However, by the concentration of measure phenomenon, we expect that for large $N$
this relation between the eigenspaces concentrates on \emph{typical} or \emph{generic
positions}, and that then the asymptotic eigenvalue distribution of $p(X_N,Y_N)$ depends
in a deterministic way only on the asymptotic eigenvalue distribution of $X_N$ and on the
asymptotic eigenvalue distribution of $Y_N$. Free probability theory replaces this vague
notion of \emph{generic position} by the mathematical precise concept of \emph{freeness}
and provides general tools for calculating the asymptotic distribution of $p(X_N,Y_N)$
out of the asymptotic distribution of $X_N$ and the asymptotic distribution of $Y_N$.

\begin{figure}[h]
\includegraphics[width=2.6in]{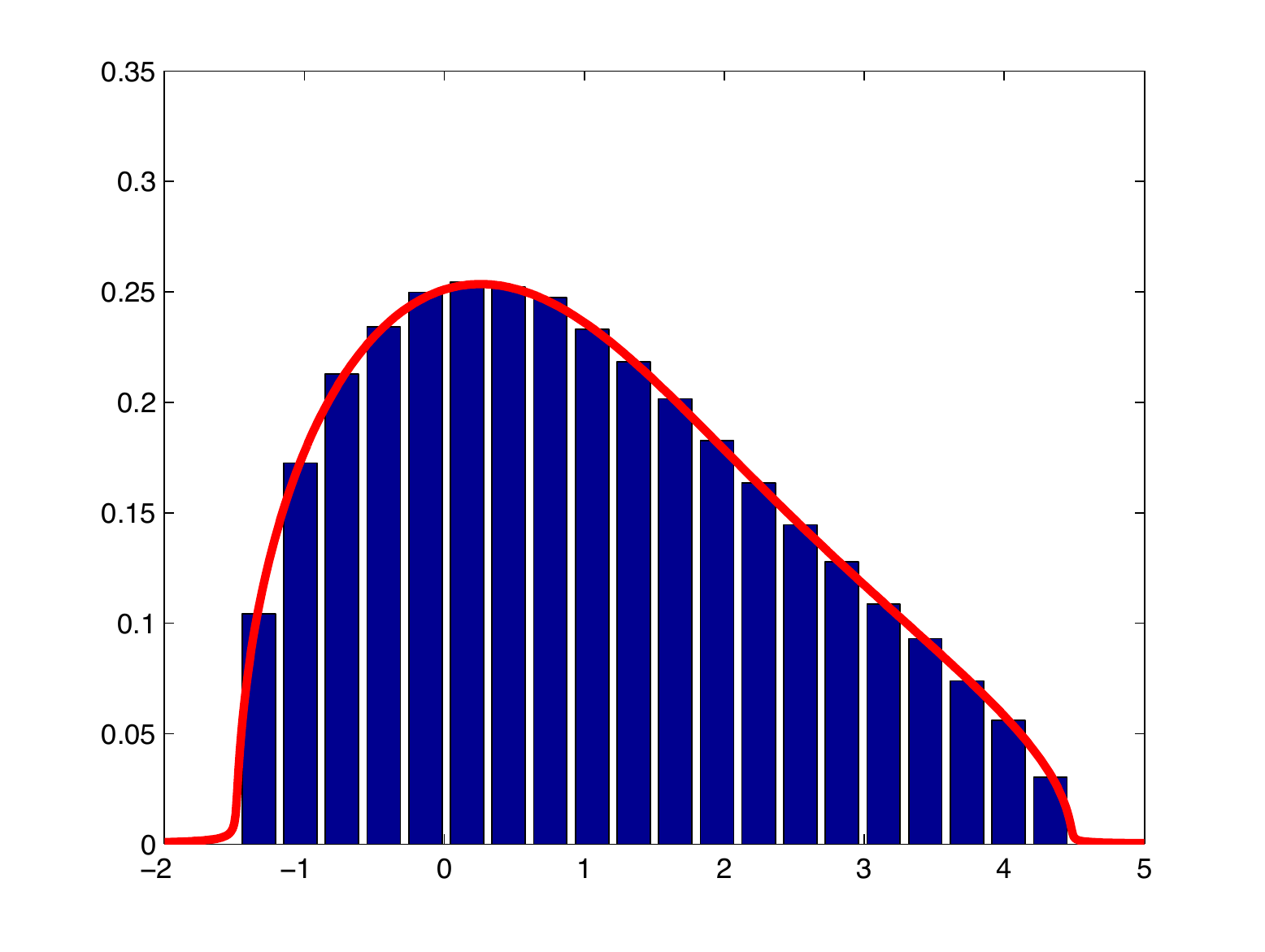}
\includegraphics[width=2.6in]{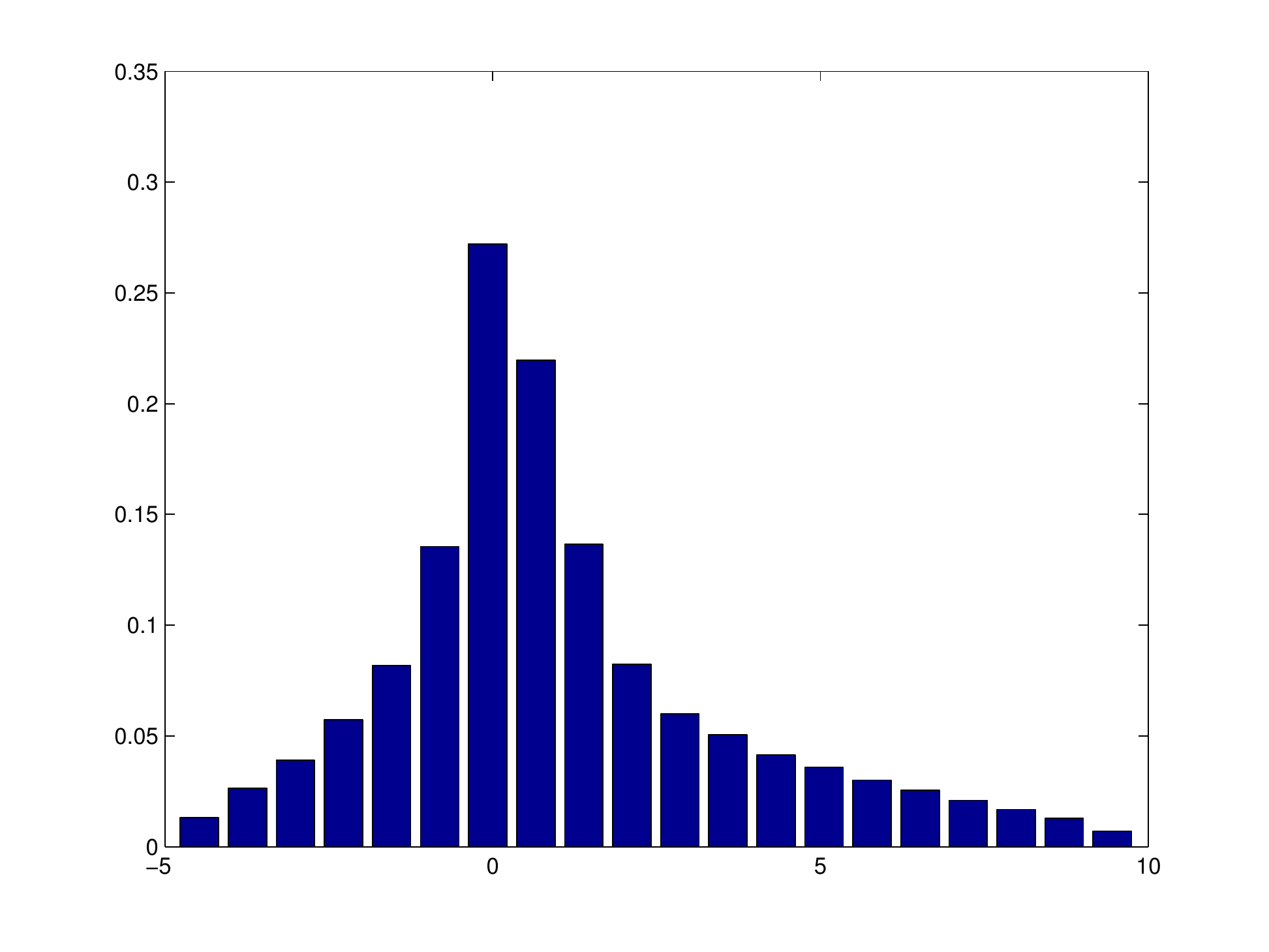}
\caption{\label{fig:sum} Histogram for a generic realization of a $3000\times 3000$ random matrix $p(X,Y)$, where $X$ and $Y$ are independent Gaussian and, respectively, Wishart random matrices: $p(X,Y)=X+Y$ (left); 
$p(X,Y)=XY+YX+X^2$ (right). In the left case, the asymptotic eigenvalue distribution is relatively easy to calculate; in the right case, no such solution was known, this case will be reconsidered in Fig.~\ref{fig:final}.} 
\end{figure}

One can convince oneself easily of the
almost sure convergence to a deterministic limit by simulations. Two examples are shown in Fig.~\ref{fig:sum}.
Actually, usually it is also not too hard to prove this almost sure convergence by appropriate variance estimates. However, what is not clear at all is the
calculation of the form of this limit shape in general. In some cases, like the left example of Fig.~\ref{fig:sum}, this was known, but in others, like the right example of Fig.~\ref{fig:sum}, no solution was known.

Our goal is to get a \emph{conceptual way of understanding} the
asymptotic eigenvalue distributions in general and also to find an \emph{algorithm for calculating} the corresponding asymptotic eigenvalue distributions.

\subsection{The moment method}

There are different methods to analyze limit distributions of random matrices. One technique, analytical in nature, is the so called resolvent method.
The main idea of this method is to derive an equation for the resolvent of the limit distribution. The advantage of this method is that there is a powerful complex analysis machinery to deal with such equations. This method also allows to look at eigenvalue distributions without finite moments. Its drawback, however, is that one cannot deal uniformly with all polynomials in $X$ and $Y$; one has to treat each $p(X,Y)$ separately. On the other side, there is a more combinatorical technique, the so-called \emph{moment method}, for calculating the limit distribution. This has the advantage that it allows, in principle, to
deal in a uniform way with all polynomials in $X$ and $Y$. In the following we will first concentrate on the moment method
in order to motivate the concept of freeness. Later we will then come back to
more analytic questions.

By $\tr (A)$ we denote the normalized trace of an $N\times N$ matrix $A$. If we want to understand the eigenvalue distribution of a selfadjoint matrix $A$, it suffices
to know the trace 
$\tr(A^k)$ of all powers of $A$: because of the invariance of the trace under conjugation with unitaries, we have for $k\in\NN$ that
$\frac 1N\bigl(\lambda_1^k+\cdots+\lambda_N^k\bigr)
= \tr(A^k)$, where $\lambda_i$ are the eigenvalues of $A$.
Therefore, instead of studying the eigenvalue distribution of a matrix $A$ directly, the moment method looks at traces of powers, $\tr(A^k)$.

Consider now our sequences of random matrices $X_N$ and $Y_N$, each of which is assumed to have almost surely an
asymptotic eigenvalue distribution. We want to understand, in the limit $N\to\infty$, the
eigenvalue distribution of $p(X_N,Y_N)$, not just for one $p$, but for all non-commutative polynomials. By the moment method, this asks for the investigation of the limit
$N\to\infty$ of $\tr\bigl(p(X_N,Y_N)^k\bigr)$ for all $k\in\NN$ and all polynomials $p$. Then it is clear that the basic objects which we have to
understand in this approach are the asymptotic \emph{mixed moments}
\begin{equation}\label{eq:mixed-moments}
\lim_{N\to\infty}\tr(X_N^{n_1}Y_N^{m_1}\cdots X_N^{n_k}Y_N^{m_k})\qquad (k\in\NN;\,
n_1,\dots,n_k,m_1,\dots,m_k\in\NN).
\end{equation}

Thus our fundamental problem is the following. If $X_N$ and $Y_N$ each have an asymptotic
eigenvalue distribution, and if $X_N$ and $Y_N$ are in generic position, do the
asymptotic mixed moments \eqref{eq:mixed-moments} 
exist? If so, can we express them in a deterministic way in terms of
the \emph{individual moments}
\begin{equation*}
\left(\lim_{N\to\infty}\tr(X_N^k)\right)_{k\in\NN}\qquad\text{and}\qquad
\left(\lim_{N\to\infty}\tr(Y_N^k)\right)_{k\in\NN}?
\end{equation*}
In order to get an idea how this might look like in a generic situation, we will consider the simplest case of two such random matrices.

\subsection{The example of two independent Gaussian random matrices}

Consider, for example, $N\times N$ random matrices $X_N$ and $Y_N$ such that
$X_N$ and $Y_N$ have asymptotic eigenvalue distributions for $N\to\infty$,
$X_N$ and $Y_N$ are independent
(i.e., the entries of $X_N$ are independent from the entries of $Y_N$) and
$Y_N$ is an unitarily invariant ensemble
(i.e., the joint distribution of its entries does not change under unitary conjugation and thus, for any unitary $N\times
N$-matrix $U_N$, $U_NY_NU_N^*$ is equivalent to the original ensemble $Y_N$ in all relevant aspects). But then we can use this $U_N$ to
rotate the eigenspaces of $Y_N$ against those of $X_N$ into a generic position, thus for
typical realizations of $X_N$ and $Y_N$ the eigenspaces should be in a generic position.

The simplest example of two such random matrix ensembles are two independent Gaussian
random matrices $X_N$ and $Y_N$. In this case one can calculate everything concretely: in
the limit $N\to\infty$, $\tr(X_N^{n_1}Y_N^{m_1}\cdots X_N^{n_k}Y_N^{m_k})$ is almost
surely given by the number of non-crossing (aka planar) pairings of the word
$$\underbrace{X\cdot X\cdots X}_{\text{$n_1$-times}}\cdot
\underbrace{Y\cdot Y\cdots Y}_{\text{$m_1$-times}}\cdots \underbrace{X\cdot X\cdots
X}_{\text{$n_k$-times}}\cdot \underbrace{Y\cdot Y\cdots Y}_{\text{$m_k$-times}},
$$
in two letters $X$ and $Y$, such that no $X$ is paired with a $Y$. 
A pairing is a decomposition of the word into pairs
of letters; if we connect the two letters from each pair by a line, drawn in the
half-plane below the word, then ``non-crossing'' means that we can do this without getting
crossings between lines for different pairs.

For example, we have 
$\lim_{N\to\infty}\tr(X_NX_NY_NY_NX_NY_NY_NX_N)=2$
because there are two such non-crossing pairings: 
$$\setlength{\unitlength}{.3cm}
\text{
\begin{picture}(8,3)
\thicklines \put(0,0){\line(0,1){2}} \put(0,0){\line(1,0){1}} \put(1,0){\line(0,1){2}}
\put(2,0){\line(0,1){2}} \put(2,0){\line(1,0){1}} \put(3,0){\line(0,1){2}}
\put(4,0){\line(0,1){2}} \put(4,0){\line(1,0){3}} \put(7,0){\line(0,1){2}}
\put(5,1){\line(0,1){1}} \put(5,1){\line(1,0){1}} \put(6,1){\line(0,1){1}}
\put(0,2.7){\makebox(0,0){$X$}} \put(1,2.7){\makebox(0,0){$X$}}
\put(2,2.7){\makebox(0,0){$Y$}} \put(3,2.7){\makebox(0,0){$Y$}}
\put(4,2.7){\makebox(0,0){$X$}} \put(5,2.7){\makebox(0,0){$Y$}}
\put(6,2.7){\makebox(0,0){$Y$}} \put(7,2.7){\makebox(0,0){$X$}}
\end{picture}}\qquad\qquad
\text{
\begin{picture}(6,3)
\thicklines \put(0,0){\line(0,1){3}} \put(0,0){\line(1,0){7}} \put(7,0){\line(0,1){3}}
\put(5,1){\line(0,1){2}} \put(5,1){\line(1,0){1}} \put(6,1){\line(0,1){2}}
\put(1,1){\line(0,1){2}} \put(1,1){\line(1,0){3}} \put(4,1){\line(0,1){2}}
\put(2,2){\line(1,0){1}} \put(2,2){\line(0,1){1}} \put(3,2){\line(0,1){1}}
\put(0,3.7){\makebox(0,0){$X$}} \put(1,3.7){\makebox(0,0){$X$}}
\put(2,3.7){\makebox(0,0){$Y$}} \put(3,3.7){\makebox(0,0){$Y$}}
\put(4,3.7){\makebox(0,0){$X$}} \put(5,3.7){\makebox(0,0){$Y$}}
\put(6,3.7){\makebox(0,0){$Y$}} \put(7,3.7){\makebox(0,0){$X$}}
\end{picture}}\qquad\qquad
$$

After some contemplation, one realizes that the above combinatorial description of the limit of $\tr(X_N^{n_1}Y_N^{m_1}\cdots X_N^{n_k}Y_N^{m_k})$ implies that the trace of a
corresponding product of centered powers,
\begin{multline}\label{eq:freeness-rm}
\lim_{N\to\infty}\tr\Big(\bigl(X_N^{n_1}-\lim_{M\to\infty}\tr(X_M^{n_1})\cdot1\bigr)\cdot
\bigl(Y_N^{m_1}-\lim_{M\to\infty}\tr(Y_M^{m_1})\cdot1\bigr)\cdots
\\ \cdots \bigl(X_M^{n_k}-\lim_{M\to\infty}\tr(X_M^{n_k})\cdot1\bigr)
\cdot\bigl(Y_N^{m_k}-\lim_{M\to\infty}\tr(Y_M^{m_k})\cdot1\bigr)\Big)
\end{multline}
is given by the number of non-crossing pairings as above, but with the additional property that  each group of $X^{n_i}$ must be connected by at least one pair to some other group $X^{n_j}$ (with $i\not =j$) and, in the same way, each group of $Y$'s must be connected to some
other group of $Y$'s. However, since the groups of $X$'s and the groups of $Y$'s are alternating, it is obvious that if we want to connect the groups in this way we will necessarily get crossings between some pairs. Thus there are actually no pairings of the required form
and we have that the term \eqref{eq:freeness-rm} is equal to zero.

One might wonder what advantage is gained by trading the explicit formula for mixed
moments \eqref{eq:mixed-moments} of independent Gaussian random matrices for the implicit relations
\eqref{eq:freeness-rm}? The drawback of the explicit formula is that the asymptotic formula for
$\tr(X_N^{n_1}Y_N^{m_1}\cdots X_N^{n_k}Y_N^{m_k})$ will be different for different random
matrix ensembles (and in many cases an explicit formula fails to exist). However, the
vanishing of \eqref{eq:freeness-rm} remains valid for many matrix ensembles. The
vanishing of \eqref{eq:freeness-rm} gives a precise meaning to our idea that the random
matrices should be in generic position; it constitutes Voiculescu's definition of
asymptotic freeness.

\begin{definition} Two sequences of matrices $(X_N)_{N\in\NN}$ and
$(Y_N)_{N\in\NN}$ are \emph{asymptotically free} if we have the vanishing of
\eqref{eq:freeness-rm} for all $k\geq 1$ and all $n_1,m_1,$ $\dots$,$ n_k,m_k\geq 1$.
\end{definition}

Provided with this definition, the intuition that unitarily invariant random matrices should
give rise to generic situations becomes now a rigorous theorem. This basic observation was
proved by Voiculescu in 1991.

\begin{theorem}[Voiculescu \cite{V-inventiones}]
\label{thm:random-matrix-freeness}
Consider $N\times N$ random matrices $X_N$ and $Y_N$ such
that: both $X_N$ and $Y_N$ have almost surely an asymptotic eigenvalue distribution for
$N\to\infty$; $X_N$ and $Y_N$ are independent; $Y_N$ is a unitarily invariant ensemble.
Then, $X_N$ and $Y_N$ are almost surely asymptotically free.
\end{theorem}

Extensions of this statement to other classes of random matrices can, for example, be found in \cite{AGZ,Dyk,MSp}.

We are now ready to give, in the next chapter, a more abstract definition of freeness. We will also see how it allows to deduce mixed moments from the individual moments.

\section{Free probability and non-crossing partitions}

\subsection{Freeness}
The starting point of free probability was the definition of freeness, given by Voiculescu in $1983$. However, this happened in the context of operator algebras, related to the isomorphism problem of free group factors. A few years later, in 1991, Voiculescu discovered the relation between random matrices and free probability, as outlined in the last chapter. These connections between operator algebras and random matrices 
led, among others, to deep results on free group factors.
In 1994, I developped a combinatorial theory of freeness, based on free cumulants; many consequences of this approach were worked out together with
Nica, see \cite{Nic, NS}.
In the following we concentrate on this combinatorial way of understanding freeness.

\begin{definition}
A pair $(\cA,\ff)$ consisting of a unital algebra $\cA$ and a linear functional
$\ff:\cA\to\CC$ with $\ff(1)=1$ is called a \emph{non-commutative probability space}.
Often the adjective ``non-commutative'' is just dropped. Elements from $\cA$ are
addressed as \emph{(non-commutative) random variables}, the numbers $\ff(a_{i(1)}\cdots
a_{i(n)})$ for such random variables $a_1,\dots,a_k\in\cA$ are called \emph{moments}, the
collection of all moments is called the \emph{joint distribution of $a_1,\dots,a_k$}.
\end{definition}

\begin{definition}
Let $(\cA,\ff)$ be a non-commutative probability space and let $I$ be an index set.

1)  Let, for each $i\in I$, $\cA_i\subset \cA$, be a unital subalgebra.
The subalgebras $(\cA_i)_{i\in I}$ are called \emph{free} or \emph{freely independent}, 
if $\ff(a_1\cdots
a_k)=0$ whenever we have: $k$ is a positive integer; $a_j\in\cA_{i(j)}$ (with $i(j)\in
I$) for all $j=1,\dots,k$; $\ff(a_j)=0$ for all $j=1,\dots,k$; and neighboring elements
are from different subalgebras, i.e., $i(1)\not=i(2), i(2)\not= i(3),\dots,
i(k-1)\not=i(k)$.

2) Let, for each $i\in I$, $x_i\in\cA$. The random variables $(x_i)_{i\in I}$ are called
\emph{free} or \emph{freely independent}, if their generated unital subalgebras are free,
i.e., if $(\cA_i)_{i\in I}$ are free, where, for each $i\in I$, $\cA_i$ is the unital
subalgebra of $\cA$ which is generated by $x_i$.
In the same way, subsets $(\cX_i)_{i\in I}$ of $\cA$ are free, if their generated unital subalgebras are so.
\end{definition}

Freeness between $x$ and $y$ is, by definition, an infinite set of equations relating various moments in
$x$ and $y$. However, one should notice
that freeness between $x$ and $y$ is actually a 
\emph{rule for calculating mixed moments} in $x$ and $y$ from
the moments of $x$ and the moments of $y$.
In this sense, freeness is
analogous to the concept of independence for classical random variables.
Hence freeness is also called \emph{free independence}.
Free probability theory investigates
these freeness relations abstractly, inspired by the philosophy that freeness should be
considered and treated as a non-commutative analogue of the classical notion of
independence.

The following examples show some calculations of mixed moments. That this works also
in general should be clear.

\begin{example}
Let us calculate, for $m,n\geq 1$, the mixed moment $\ff (x^ny^m)$ of some free random variables $x$ and $y$.
By the definition of freeness it follows that 
$\ff[(x^n-\ff(x^n)1)(y^m-\ff(y^m)1)]=0$. This gives
$$\ff(x^ny^m)-\ff(x^n\cdot 1)\ff(y^m)-\ff(x^n)\ff(1\cdot
y^m)+\ff(x^n)\ff(y^m)\ff(1\cdot 1)=0,$$
and hence
$\ff(x^ny^m)=\ff(x^n)\cdot \ff(y^m)$.
\end{example}

The above is the same result as for independent classical random variables. However, this is misleading.
Free independence is a different rule from classical independence; free
independence occurs typically for \emph{non-commuting
random variables}, like operators on Hilbert spaces or (random) matrices.

\begin{example}
Let $x$ and $y$ be some free random variables. By definition of freeness we get
$$\ff[(x-\ff(x)1)\cdot (y-\ff(y)1)\cdot
(x-\ff(x)1)\cdot (y-\ff(y)1)]=0,$$ 
which results after some elementary, but lengthy calculations and many cancellations in
\begin{equation}\label{eq:xyxy}
\ff(xyxy)=\ff(xx)\cdot\ff(y)\cdot\ff(y)+\ff(x)\cdot\ff(x)\cdot\ff(yy)
-\ff(x)\cdot\ff(y)\cdot\ff(x)\cdot \ff(y).
\end{equation}
\end{example}

We see that this result is different from the one for independent classical
(and thus commuting) random variables. 
It is important to note that freeness plays a
similar role in the non-commutative world as independence plays in the classical world, but
that freeness is not a generalization of independence: independent random variables can be free
only in very trivial situations. Freeness is a theory for genuinely non-commuting random variables.

\subsection{Understanding the freeness rule: the idea of cumulants}
The main idea in this section is to write moments in terms of other quantities, which we call \emph{free cumulants}.
We will see that freeness is much easier to describe on the level of free cumulants, namely by the vanishing of mixed cumulants. There is also a nice relation between moments and cumulants, given by summing over
\emph{non-crossing or planar partitions}.

\begin{definition}
1) A \emph{partition} of $\{1,\dots,n\}$ is a decomposition
$\pi=\{V_1,\dots,V_r\}$ of $\{1,\dots,n\}$ into subsets $V_i$ 
with 
$$V_i\not=\emptyset,\qquad V_i\cap V_j=\emptyset\quad (i\not=j),\qquad
\bigcup_i V_i=\{1,\dots,n\}.$$
The $V_i$ are called the \emph{blocks} of $\pi$.
The set of all such partitions is denoted by $\cP(n)$.

2) A partition $\pi$ is \emph{non-crossing} if we do not have
$p_1<q_1<p_2<q_2$
such that $p_1,p_2$ are in a same block, $q_1,q_2$ are in a same block, but those two blocks are different.
By $NC(n)$ we will denote the set of all non-crossing partions of $\{1,\dots,n\}$.
\end{definition}

Let us remark that
$NC(n)$ is actually a lattice with respect to refinement order.

\begin{definition}
For a unital linear functional
$\ff:\cA\to\CC$
we define the \emph{free cumulants} $\kk_n:\cA^n\to\CC$ (for all
$n\geq 1$)
as multi-linear functionals
by the \emph{moment-cumulant relation}
$$\ff(a_1\cdots a_n)=\sum_{\pi\in NC(n)}\kk_\pi(a_1,\dots, a_n).$$
Here, $\kk_\pi$ is a product of cumulants: one term for each block of $\pi$,
and the arguments are given by the elements corresponding to the respective blocks. This, as well as the fact that these equations define the free cumulants uniquely, will be illustrated by the following examples.
\end{definition}

This definition is motivated by a similar formula for classical cumulants. The only difference is that in the classical case non-crossing partitions $NC(n)$ are replaced by all partitions $\cP(n)$.

\setlength{\unitlength}{0.2cm}
  \newsavebox{\NCianew}
  \savebox{\NCianew}(0,1){
  \thicklines
  \put(0,0){\line(0,1){1}}}
\setlength{\unitlength}{0.2cm}
  \newsavebox{\NCiianew}
  \savebox{\NCiianew}(1,1){
  \thicklines
  \put(0,0){\line(0,1){1}}
  \put(0,0){\line(1,0){1}}
  \put(1,0){\line(0,1){1}}}
  \newsavebox{\NCiibnew}
  \savebox{\NCiibnew}(1,1){
  \thicklines
  \put(0,0){\line(0,1){1}}
  \put(1,0){\line(0,1){1}}}
\setlength{\unitlength}{0.2cm}
  \newsavebox{\NCiibnewa}
  \savebox{\NCiibnewa}(1,1){
  \thicklines
  \put(0,0){\line(0,1){1}}
  \put(1,0){\line(0,1){1}}}

  \setlength{\unitlength}{2cm}
  \newsavebox{\NCiiiannewa}
  \savebox{\NCiiiannewa}(2,2){
  \thicklines
  \put(0,0){\line(0,1){2}}
  \put(0,0){\line(1,0){2}}
  \put(1,0){\line(0,1){2}}
  \put(2,0){\line(0,1){2}}}
  \newsavebox{\NCiiibnnewa}
  \savebox{\NCiiibnnewa}(2,2){
  \thicklines
  \put(0,0){\line(0,1){2}}
  \put(1,0){\line(0,1){2}}
  \put(1,0){\line(1,0){1}}
  \put(2,0){\line(0,1){2}}}
  \newsavebox{\NCiiicnnewa}
  \savebox{\NCiiicnnewa}(2,2){
  \thicklines
  \put(0,0){\line(0,1){2}}
  \put(0,0){\line(1,0){1}}
  \put(1,0){\line(0,1){2}}
  \put(2,0){\line(0,1){2}}}
  \newsavebox{\NCiiidnnewa}
  \savebox{\NCiiidnnewa}(2,2){
  \thicklines
  \put(0,0){\line(0,1){2}}
  \put(0,0){\line(1,0){2}}
  \put(2,0){\line(0,1){2}}
  \put(1,1){\line(0,1){1}}}
  \newsavebox{\NCiiiennewa}
  \savebox{\NCiiiennewa}(2,2){
  \thicklines
  \put(0,0){\line(0,1){2}}
  \put(1,0){\line(0,1){2}}
  \put(2,0){\line(0,1){2}}}

\begin{example}
Let us calculate some examples for cumulants for small $n$.

For $n=1$ there exists only one partition, $\usebox{\NCianew}$ , so that the first moment and the first cumulant are the same: $\ff(a_1) = \kk_1(a_1)$.

For $n=2$ there are two partitions,
$\usebox{\NCiianew}\quad \text{and} \usebox{\NCiibnew}$ ,
and both are non-crossing. By the moment-cumulant formula we get 
$\ff(a_1a_2)=\kk_2(a_1,a_2)+\kk_1(a_1)\kk_1(a_2)$,
and thus $\kk_2$ is nothing but the covariance
$\kk_2(a_1,a_2)=\ff(a_1a_2)-\ff(a_1)\ff(a_2)$.

\setlength{\unitlength}{0.25cm}
  \newsavebox{\NCiiiannew}
  \savebox{\NCiiiannew}(2,2){
  \thicklines
  \put(0,0){\line(0,1){2}}
  \put(0,0){\line(1,0){2}}
  \put(1,0){\line(0,1){2}}
  \put(2,0){\line(0,1){2}}}
  \newsavebox{\NCiiibnnew}
  \savebox{\NCiiibnnew}(2,2){
  \thicklines
  \put(0,0){\line(0,1){2}}
  \put(1,0){\line(0,1){2}}
  \put(1,0){\line(1,0){1}}
  \put(2,0){\line(0,1){2}}}
  \newsavebox{\NCiiicnnew}
  \savebox{\NCiiicnnew}(2,2){
  \thicklines
  \put(0,0){\line(0,1){2}}
  \put(0,0){\line(1,0){1}}
  \put(1,0){\line(0,1){2}}
  \put(2,0){\line(0,1){2}}}
  \newsavebox{\NCiiidnnew}
  \savebox{\NCiiidnnew}(2,2){
  \thicklines
  \put(0,0){\line(0,1){2}}
  \put(0,0){\line(1,0){2}}
  \put(2,0){\line(0,1){2}}
  \put(1,1){\line(0,1){1}}}
  \newsavebox{\NCiiiennew}
  \savebox{\NCiiiennew}(2,2){
  \thicklines
  \put(0,0){\line(0,1){2}}
  \put(1,0){\line(0,1){2}}
  \put(2,0){\line(0,1){2}}}
	
In the same recursive way, we are able to compute the third cumulant. There are five partitions of the set of three elements:
$$\usebox{\NCiiiannew} \qquad\usebox{\NCiiibnnew}\qquad \usebox{\NCiiicnnew}\qquad
               \usebox{\NCiiidnnew}\qquad \usebox{\NCiiiennew}$$
Still, they are all non-crossing and the moment-cumulant formula gives
\begin{align*}
\ff(a_1a_2a_3)&=
\kk_3(a_1,a_2,a_3)
+\kk_1(a_1)\kk_2(a_2,a_3)\\&\quad+\kk_2(a_1,a_2)\kk_1(a_3)
+\kk_2(a_1,a_3)\kk_1(a_2)  
+\kk_1(a_1)\kk_1(a_2)\kk_1(a_3)
\end{align*}
and hence
\begin{multline*}
\kk_3(a_1,a_2,a_3)=\ff(a_1a_2a_3)-\ff(a_1)\ff(a_2a_3)-\ff(a_1a_2)\ff(a_3)\\
-\ff(a_1a_3)\ff(a_2)+2 \ff(a_1)\ff(a_2)\ff(a_3).
\end{multline*}

\setlength{\unitlength}{0.25cm}
\newsavebox{\NCa}
\savebox{\NCa}(3,2){ \thicklines \put(0,0){\line(0,1){2}} \put(0,0){\line(1,0){3}}
\put(1,0){\line(0,1){2}} \put(2,0){\line(0,1){2}} \put(3,0){\line(0,1){2}}}

\newsavebox{\NCb}
\savebox{\NCb}(3,2){ \thicklines \put(0,0){\line(0,1){2}} \put(1,0){\line(0,1){2}}
\put(1,0){\line(1,0){2}} \put(2,0){\line(0,1){2}} \put(3,0){\line(0,1){2}}}

\newsavebox{\NCc}
\savebox{\NCc}(3,2){ \thicklines \put(0,0){\line(0,1){2}} \put(0,0){\line(1,0){3}}
\put(1,1){\line(0,1){1}} \put(2,0){\line(0,1){2}} \put(3,0){\line(0,1){2}}}

\newsavebox{\NCd}
\savebox{\NCd}(3,2){ \thicklines \put(0,0){\line(0,1){2}} \put(0,0){\line(1,0){3}}
\put(1,0){\line(0,1){2}} \put(2,1){\line(0,1){1}} \put(3,0){\line(0,1){2}}}

\newsavebox{\NCe}
\savebox{\NCe}(3,2){ \thicklines \put(0,0){\line(0,1){2}} \put(0,0){\line(1,0){2}}
\put(1,0){\line(0,1){2}} \put(2,0){\line(0,1){2}} \put(3,0){\line(0,1){2}}}

\newsavebox{\NCf}
\savebox{\NCf}(3,2){ \thicklines \put(0,0){\line(0,1){2}} \put(0,0){\line(1,0){1}}
\put(1,0){\line(0,1){2}} \put(2,0){\line(0,1){2}} \put(3,0){\line(0,1){2}}
\put(2,0){\line(1,0){1}}}

\newsavebox{\NCg}
\savebox{\NCg}(3,2){ \thicklines \put(0,0){\line(0,1){2}} \put(0,0){\line(1,0){3}}
\put(1,1){\line(0,1){1}} \put(2,1){\line(0,1){1}} \put(3,0){\line(0,1){2}}
\put(1,1){\line(1,0){1}}}

\newsavebox{\NCh}
\savebox{\NCh}(3,2){ \thicklines \put(0,0){\line(0,1){2}} \put(1,0){\line(0,1){2}}
\put(2,0){\line(0,1){2}} \put(2,0){\line(1,0){1}} \put(3,0){\line(0,1){2}}}

\newsavebox{\NCi}
\savebox{\NCi}(3,2){ \thicklines \put(0,0){\line(0,1){2}} \put(1,0){\line(0,1){2}}
\put(2,0){\line(0,1){2}} \put(1,0){\line(1,0){1}} \put(3,0){\line(0,1){2}}}

\newsavebox{\NCj}
\savebox{\NCj}(3,2){ \thicklines \put(0,0){\line(0,1){2}} \put(1,0){\line(0,1){2}}
\put(2,0){\line(0,1){2}} \put(0,0){\line(1,0){1}} \put(3,0){\line(0,1){2}}}

\newsavebox{\NCk}
\savebox{\NCk}(3,2){ \thicklines \put(0,0){\line(0,1){2}} \put(1,0){\line(0,1){2}}
\put(2,1){\line(0,1){1}} \put(1,0){\line(1,0){2}} \put(3,0){\line(0,1){2}}}

\newsavebox{\NCl}
\savebox{\NCl}(3,2){ \thicklines \put(0,0){\line(0,1){2}} \put(1,1){\line(0,1){1}}
\put(2,1){\line(0,1){1}} \put(0,0){\line(1,0){3}} \put(3,0){\line(0,1){2}}}

\newsavebox{\NCm}
\savebox{\NCm}(3,2){ \thicklines \put(0,0){\line(0,1){2}} \put(1,1){\line(0,1){1}}
\put(2,0){\line(0,1){2}} \put(0,0){\line(1,0){2}} \put(3,0){\line(0,1){2}}}

\newsavebox{\NCn}
\savebox{\NCn}(3,2){ \thicklines \put(0,0){\line(0,1){2}} \put(1,0){\line(0,1){2}}
\put(2,0){\line(0,1){2}} \put(3,0){\line(0,1){2}}}

The first difference to the classical theory occurs now for $n=4$; there are 15 partitions of the set of four elements, but one is crossing and there are only 14 non-crossing partitions:
$$\usebox{\NCa}\qquad
\usebox{\NCb}\qquad\usebox{\NCc}\qquad\usebox{\NCd}
\qquad\usebox{\NCe}
\qquad\usebox{\NCf}\qquad\usebox{\NCg}
$$
$$
\usebox{\NCh}\qquad
\usebox{\NCi}
\qquad
\usebox{\NCj}
\qquad
\usebox{\NCk}\qquad
\usebox{\NCl}\qquad
\usebox{\NCm}\qquad
\usebox{\NCn}
$$
Hence the moment-cumulant formula yields
\begin{align*}
\ff(a_1a_2&a_3a_4)=
\kk_4(a_1,a_2,a_3,a_4)+ \kk_1(a_1)\kk_3(a_2,a_3,a_4)+
\kk_1(a_2)\kk_3(a_1,a_3,a_4)\\&+\kk_1(a_3)\kk_3(a_1,a_2,a_4)+
\kk_3(a_1,a_2,a_3)\kk_1(a_4)
+\kk_2(a_1,a_2)\kk_2(a_3,a_4)\\&+\kk_2(a_1,a_4)\kk_2(a_2,a_3)
+\kk_1(a_1)\kk_1(a_2)\kk_2(a_3,a_4)
+\kk_1(a_1)\kk_2(a_2,a_3)\kk_1(a_4)\\&+\kk_2(a_1,a_2)\kk_1(a_3)\kk_1(a_4)+\kk_1(a_1)
\kk_2(a_2,a_4)\kk_1(a_3)
+\kk_2(a_1,a_4)\kk_1(a_2)\kk_1(a_3)\\&
+\kk_2(a_1,a_3)\kk_1(a_2)\kk_1(a_4)+\kk_1(a_1)\kk_1(a_2)
\kk_1(a_3)\kk_1(a_4).
\end{align*}
As before, this can be resolved for $\kk_4$ in terms of moments. 
\end{example}

The reader should by now be convinced that one can actually rewrite the moment-cumulant equations also the other way round as cumulant-moment equations. More precisely, this can be achieved by a M\"obius inversion on the poset of non-crossing partitions resulting in
$$\kk_n(a_1,\dots,a_n)=\sum_{\pi\in NC(n)}\ff_\pi(a_1,\dots,a_n)
\mu(\pi,1_n),$$
where $\ff_\pi$ is a product of moments according to the block structure of
$\pi$ and $\mu$ is the M\"obius function of $NC(n)$.

Whereas $\kk_1$,
$\kk_2$, and $\kk_3$ are the same as the corresponding classical cumulants, the free cumulant $\kk_4$ and all the higher ones are different from their classical counterparts.

\subsection{Freeness corresponds to vanishing of mixed
cumulants}

The following theorem shows that freeness is much easier to describe on the level of cumulants than on the level of moments. This characterization is at the basis of most calculations with free cumulants.

\begin{theorem}[Speicher \cite{Sp-mult}]\label{thm:vanishing}
The fact that $x$ and $y$ are free is equivalent to the fact that
$\kk_n(a_1,\dots,a_n)=0$
whenever:
$n\geq 2$,
$a_i\in\{x,y\}$ for all $i$, and
there are at least two indices $i,j$ such that $a_i=x$ and $a_j=y$.
\end{theorem}

A corresponding statement is also true for more than two random variables: freeness is equivalent to the vanishing of mixed cumulants.

\begin{example}
If $x$ and $y$ are free, then we have
$$\ff(xyxy)=\kk_1(x)\kk_1(x)\kk_2(y,y)+
\kk_2(x,x)\kk_1(y)\kk_1(y)+\kk_1(x)\kk_1(y)\kk_1(x)\kk_1(y),$$
corresponding to the three non-crossing partitions of $xyxy$ which connect $x$ only with $x$ and $y$ only with $y$:  \setlength{\unitlength}{.25cm}
$$\qquad\text{
\begin{picture}(8,3)
\thicklines {\put(0,0){\line(0,1){2}}
\put(2,1){\line(0,1){1}}}{ \put(1,0){\line(1,0){1}}
\put(1,0){\line(0,1){2}}
 \put(2,0){\line(1,0){1}} \put(3,0){\line(0,1){2}}}
\put(0,2.7){\makebox(0,0){$x$}} \put(1,2.7){\makebox(0,0){$y$}}
\put(2,2.7){\makebox(0,0){$x$}} \put(3,2.7){\makebox(0,0){$y$}} \end{picture}}
\qquad\text{
\begin{picture}(8,3)
\thicklines  {\put(0,0){\line(1,0){2}}
\put(0,0){\line(0,1){2}}
\put(2,0){\line(0,1){2}}}{\put(3,0){\line(0,1){2}}
\put(1,1){\line(0,1){1}}} \put(0,2.7){\makebox(0,0){$x$}}
\put(1,2.7){\makebox(0,0){$y$}} \put(2,2.7){\makebox(0,0){$x$}}
\put(3,2.7){\makebox(0,0){$y$}}
\end{picture}}\qquad\qquad \text{
\begin{picture}(8,3)
\thicklines  {\put(0,0){\line(0,1){2}}
\put(2,0){\line(0,1){2}}}{\put(1,0){\line(0,1){2}}
\put(3,0){\line(0,1){2}}} \put(0,2.7){\makebox(0,0){$x$}}
\put(1,2.7){\makebox(0,0){$y$}} \put(2,2.7){\makebox(0,0){$x$}}
\put(3,2.7){\makebox(0,0){$y$}} \end{picture}}
$$
Rewriting the cumulants in terms of moments recovers of course 
the formula \eqref{eq:xyxy}.
\end{example}

This description of freeness in terms of free cumulants is related to the planar
approximations in random matrix theory. In a sense some aspects of this theory of
freeness were anticipated (but mostly neglected) in the physics community in the paper
\cite{Cvi82}.

\section{Sum of free variables: free convolution}
Let $x$, $y$ be two free random variables.
Then, by freeness, the moments of $x+y$ are uniquely determined by
the moments of $x$ and the moments of $y$. But is there an effective way to calculate the distribution of $x+y$ if we know 
the distribution of $x$ and the distribution of $y$?

\subsection{Free convolution}
We usually consider this question in a context where we have some more analytic structure. Formally, a good frame for this is a
\emph{$C^*$-probability space} $(\cA,\ff)$, where $\cA$ is a $C^*$-algebra (i.e., a norm-closed
$*$-subalgebra of the algebra of bounded operators on a Hilbert space) and $\ff$ is a
state, i.e. it is positive in the sense $\ff(aa^*)\geq 0$ for all $a\in \cA$. Concretely
this means that our random variables can be realized as bounded operators on a Hilbert
space and $\ff$ can be written as a vector state $\ff(a)=\langle a\xi,\xi\rangle$ for
some unit vector $\xi$ in the Hilbert space.

In such a situation the distribution of a selfadjoint random variable $x$ can be identified
with a compactly supported probability measure $\mu_x$ on $\RR$, via
\begin{equation}\label{eq:mom-prob}
\ff(x^n)=\int_\RR t^nd\mu_x(t)\qquad\text{for all $n\in\NN$}.
\end{equation}

Then we say that the distribution of $x+y$, if $x$ and $y$ are free, is the
\emph{free convolution} of the
distribution $\mu _x$ of $x$ and the distribution $\mu _y$ of $y$ and denote this by
$\mu_{x+y}=\mu_x\boxplus \mu_y$.
By considering unbounded selfadjoint operators (and replacing moments of $x$ by bounded functions of $x$ in \eqref{eq:mom-prob}) one can extend this free convolution also to a binary operation on arbitrary probability measures on $\RR$, see \cite{BV}.

In principle, freeness determines $\mu_x\boxplus\mu_y$ in terms of $\mu_x$ and $\mu_y$, but the concrete nature of this connection on the level
of moments is not apriori clear.
However, by Theorem~\ref{thm:vanishing}, there is an easy rule on the level of free cumulants:  if $x$ and $y$ are free then we have for all $n\geq 1$ that
$\kk_n(x+y,x+y,\dots,x+y)=\kk_n(x,x,\dots,x)+\kk_n(y,y,\dots,y)$,
because all mixed cumulants in $x$ and $y$ vanish. 

Thus, the description of the free convolution has now been shifted to understanding the relation between moments and cumulants. A main step for this understanding is the fact that the combinatorial relation between moments and
cumulants can also be rewritten easily as a relation between corresponding formal power
series.

\subsection{Relation between moments and free cumulants}
We denote the $n$-th moment of $x$ by $m_n:=\ff (x^n)$
and the $n$-th free cumulant of $x$ by
$\kk_n:=\kk_n(x,x,\dots,x)$.
Then, the combinatorical relation between them is given by the moment-cumulant formula
\begin{equation}\label{eq:mom-cum}
m_n=\sum_{\pi\in NC(n)}\kk_\pi,
\end{equation}
where $\kk_\pi=\kk_{\vert V_1\vert}\cdots \kk_{\vert V_s\vert}$ for
$\pi=\{V_1,\dots,V_s\}$.
The next theorem shows that this combinatorial relation can be rewritten into a functional relation between the corresponding formal power series.

\begin{theorem}[Speicher \cite{Sp-mult}]\label{thm:functional}
Consider formal power series $M(z)=1+\sum_{n=1}^\infty m_n z^n$ and $C(z)=1+\sum_{n=1}^\infty \kk_n z^n$.
Then the relation 
\eqref{eq:mom-cum}
between the coefficients is equivalent to the relation
$M(z)=C[zM(z)]$.
\end{theorem}

The main step in the proof of this is to observe that a non-crossing partition can be described by its first block (i.e., the block containing the point 1) and by the non-crossing partitions of the points between the legs of the first block. This leads to 
the following recursive relation between free cumulants and moments:
$$
m_n=\sum_{s=1}^n \sum_{i_1,\dots,i_s\geq 0\atop i_1+\cdots +i_s+s=n}
\kk_s m_{i_1}\cdots m_{i_s}.$$

An early instance of the functional relation in Theorem \ref{thm:functional} appeared also in the work of Beissinger \cite{Beissinger}, for the special problem of counting non-crossing partitions by decomposing them into irreducible components.

\begin{remark}
Classical cumulants $c_k$ are combinatorially defined by the analogous
formula
$m_n=\sum_{\pi\in\cP(n)} c_\pi$.
In terms of exponential generating power series
$\tilde M(z)=1+\sum_{n=1}^\infty \frac {m_n} {n!} z^n$ and
$\tilde C(z)=\sum_{n=1}^\infty \frac {c_n}{n!} z^n$
this is equivalent to 
$\tilde C(z)=\log \tilde M(z)$.
\end{remark}

\subsection{The Cauchy transform}
For a selfadjoint random variable $x$, with corresponding probability measure $\mu_x$ according to Eq.~\eqref{eq:mom-prob}, we define the Cauchy transform $G$ by
$$G(z):=\ff\left(\frac 1{z-x}\right)=\int_\RR\frac 1{z-t}d\mu_x(t).$$
If $\mu_x$ is compactly supported we can expand this into a formal power series:
$$G(z)=\sum_{n=0}^\infty \frac{\ff(x^n)}{z^{n+1}}=\frac {M\left(\frac 1z\right)}z.$$
Therefore, on a formal level $M(z)$ and $G(z)$ contain the same information. However, $G(z)$ has many advantages over $M(z)$.
Namely, the Cauchy transform is an analytic function $G:\CC^+\to \CC^-$
and we can recover $\mu _x$ from $G$ by using the 
\emph{Stieltjes inversion formula}:
$$d\mu_x(t)=-\frac 1\pi \lim_{\ee\to 0}\Im G(t+i\ee)dt.$$
Here, $\Im$ denotes the imaginary part and the convergence in this equation is weak convergence of probability measures;
the right hand side is, for any $\ee>0$, the density of a probability measure. 

\subsection{The $R$-transform}
Voiculescu showed in \cite{DVJFA} the existence of the free cumulants of a random variable by general arguments, but without having
a combinatorial interpretation for them. There he
defined the following variant of the cumulant generating series $C(z)$. 

\begin{definition}
For a random variable $x\in\cA$ we define its \emph{$R$-transform} by 
$$R(z)=\sum_{n=1}^\infty \kk_n(x,\dots,x) z^{n-1}.$$
\end{definition}   

Then by a simple application of our last theorem we get the following result. The original 
proof of Voiculescu was much more analytical.

\begin{theorem}[Voiculescu \cite{DVJFA}, Speicher \cite{Sp-mult}]
1) For a random variable we have the relation
$\frac 1{G(z)}+R[G(z)]=z$ between its Cauchy and $R$-transform.
\\
2) If $x$ and $y$ are free, then we have $R_{x+y}(z)=R_x(z)+R_y(z)$.
\end{theorem}

\subsection{The $R$-transform as an analytic object}
In the last sections we considered the $R$-transform only as a formal power series. But for more advanced investigations as well as for explicit calculations it is necessary
to study the analytic properties of this object. It is easy to see that for bounded selfadjoint random variables
the $R$-transform can be established as an analytic function via power series expansions around the point infinity in the complex plane. But there are some problems with the analytic properties of the $R$-transform. One problem is that 
the $R$-transform can, in contrast to the Cauchy transform, in general not be
defined on all of the upper complex half-plane, but only in some truncated cones (which depend on the considered variable).
Another problem is that the equation $\frac 1{G(z)}+R[G(z)]=z$
does in general not allow explicit solutions and there is no good numerical algorithm for dealing with this. Therefore one is in need of other tools, which allow to compute free convolutions in a more efficient way.

\subsection{An alternative to the $R$-transform: subordination}

Let $x$ and $y$ be free. 
Put 
$w:=R_{x+y}(z)+1/z$, then
$$
G_{x+y}(w)=z=G_x[R_x(z)+1/z]=G_x[w-R_y(z)]=G_x[w-R_y[G_{x+y}(w)]].
$$
Thus, with
$\omega(z):=z-R_y[G_{x+y}(z)]]$,
we have the subordination
$G_{x+y}(z)=G_x\bigl(\omega(z)\bigr)$.
Though the above manipulations were just on a formal level, it turns out that this subordination function $\omega$ is, for selfadjoint $x$ and $y$ in a $C^*$-probability space, always a nice analytic object and amenable to robust calculation algorithms.
The subordination property has first been proved in
\cite{V3} by Voiculescu under a genericity assumption, and in full generality by Biane \cite{Biane1}.

It was noted, and for the first time explicitly
formulated in \cite{CG}, that the subordination property is equivalent to the $R$-transform approach, 
but has better analytic properties.
A particularly nice feature is that the subordination function can be recovered by fixed point arguments, as shown in \cite{BB07}.

\begin{theorem}[Belinschi, Bercovici \cite{BB07}]
Let $(\cA,\ff)$ be a $C^*$-probability space and let $x=x^*$ and $y=y^*$ in $\cA$ be free. Put 
$F(z):=\frac 1{G(z)}$.
Then there exists an analytic map
$\omega:\CC^+\to\CC^+$ (depending both on $x$ and $y$)
such that
$$F_{x+y}(z)=F_x\bigl(\omega(z)\bigr)\qquad\text{and}\qquad
G_{x+y}(z)=G_x\bigl(\omega(z)\bigr).$$
The subordination function $\omega(z)$ is given as the unique fixed point in the upper half-plane of the map $f_z:\CC^+\to\CC^+$, given by
$$f_z(w)=F_y(F_x(w)-w+z)-(F_x(w)-w).$$
\end{theorem}

\section{Polynomials in several random matrices}

Our original problem was to calculate the asymptotic eigenvalue distribution of selfadjoint polynomials in several
independent random matrices in generic position. We have now a conceptual grasp on this problem by relating it to free probability theory via the 
basic result of Voiculescu which tells us that our random matrices become almost surely asymptoticially free.
This allows us to reduce our random matrix problem to the problem of polynomials in free variables:
If the random matrices $X_1,\dots,X_k$ are asymptotically freely independent, then the eigenvalue distribution of a polynomial $p(X_1,\dots,X_k)$ is asymptotically given by
the distribution of $p(x_1,\dots,x_k)$, where 
$x_1,\dots,x_k$ are freely independent variables, and 
the distribution of $x_i$ is the asymptotic distribution of $X_i$.

So now the question is: Can we calculate the distribution of polynomials in free variables?
We have seen that free convolution gives effective analytic tools for dealing with the simplest polynomial, the sum of two matrices.
By using this, we calculated for example the form of the limiting eigenvalue distribution for the sum of an independent Gaussian and Wishart matrix in the left figure of Fig.~\ref{fig:sum}.
But what can we say for more general polynomials, like the one considered in the right figure of Fig.~\ref{fig:sum}.

For this problem, both from the random matrix and the free probability side, there is a long list of contributions which provide 
solutions for special choices of the polynomial $p$. In the context of free probability, Voiculescu solved it in \cite{DVJFA} and \cite{V2} for the cases of $p(x,y)=x+y$ and $p(x,y)=xy^2x$
(corresponding to the additive and multiplicative free convolution) with the introduction of the $R$- and $S$-transform, respectively. Nica and
Speicher could give in \cite{NS-commut} a solution for the problem of the free commutator, $p(x,y)=i(xy-yx)$. 

In the random matrix context, this problem
was addressed for various polynomials -- and usually, also for specific choices of the distributions of the $X_i^{(N)}$ -- by many authors, starting with the work of Marchenko-Pastur \cite{MP}. For a more extensive list of contributions in this context we refer to the books \cite{BS,CB,Gir,TV}. 
Some of those situations were also treated by operator-valued
free probability tools, see in particular \cite{BSTV, Sh1996, SpV}.

All those investigations were specific for the considered polynomial and up to now there has not
existed a master algorithm which would work for all polynomials. 

Actually, there is no hope to calculate effectively  general
polynomials in freely independent variables within the frame of free probability theory as presented up to now.
However, there is a possible way to deal with such a general situation, by the use of a linearization trick.
This trick will be the main topic of the next chapter. 

\section{The linearization trick}

The idea of this trick is: in order to understand general polynomials in non-commuting variables, it suffices to understand matrices of \emph{linear} polynomials in those variables. 
Such linearization ideas seem to be around in many different communities. In the context of 
operator algebras, Voiculescu used such a linearization philosophy as one motivation for his
work on operator-valued free probability \cite{V1995}. A seminal concrete form is due to
Haagerup and Thorbj\o rnsen \cite{HaagerupThorbjornsen2}, who used such techniques to study the largest eigenvalue of polynomials in independent Gaussian random matrices. In 2011, based on the Schur complement, Anderson \cite{Anderson1} developped a selfadjoint version of the linearization trick, which turns out to be the right tool in our context.
We present this version of Anderson in the following.

\begin{definition}
Consider a polynomial $p$ in several non-commuting variables.
A  \emph{linearization} of $p$ is an $N\times N$ matrix (with $N\in\NN$) of the form 
$$\hat p=\begin{pmatrix} 
0&u\\ v&Q
\end{pmatrix}, $$
where:
\begin{itemize}
\item
$u, v, Q$ are matrices of appropriate sizes: $u$ is $1\times (N-1)$; $v$ is $(N-1)\times 1$; and $Q$ is $(N-1)\times (N-1)$.
\item
$Q$ is invertible and we have
$p=-uQ^{-1}v$.
\item
The entries of $\hat p$ are polynomials in the variables, each of degree $\leq 1$.
\end{itemize}
\end{definition}

A linearization is of course not uniquely determined by the above requirements.
The crucial fact is that such linearizations always exist. Furthermore, they can be chosen in such  a way that they preserve selfadjointness.

\begin{theorem}[Anderson \cite{Anderson1}]
For each $p$ there exists a linearization $\hat p$ (with an explicit algorithm
for finding those). Moreover if $p$ is selfadjoint, then this $\hat p$ is also selfadjoint.
\end{theorem}

\begin{example}
We consider the selfadjoint non-commutative polynomial $p=xy+yx+x^2$. 
Then a selfadjoint linearization of $p$ is the matrix
\begin{equation}\label{eq:phat}
\hat p=\begin{pmatrix}
0&\quad&x\quad&&\frac x2+y\\
x&\quad&0\quad&&-1\\
\frac x2+y&\quad&-1\quad&&0
\end{pmatrix},
\end{equation}
because we have
$$\begin{pmatrix}
x& \frac x2+y
\end{pmatrix}
\begin{pmatrix}
0&-1\\-1&0
\end{pmatrix}^{-1}
\begin{pmatrix}
x\\ \frac x2+y
\end{pmatrix}
=-(xy+yx+x^2).$$
\end{example}

At this point it might not be clear what this linearization trick has to do with our problem.
What we are interested in is the distribution of $p$, which can be recovered from the Cauchy transform of $p$, which is given by taking expectations of resolvents of $p$. Thus we need
control of inverses of $z-p$. How can the linearization $\hat p$ give information
on those?

For $z\in\CC$ we put 
$b=\begin{tiny}{\begin{pmatrix}
z&0\\ 0&0
\end{pmatrix}}
\end{tiny}$ and then it follows
$$b-\hat p=\begin{pmatrix} 
z&-u\\ -v&-Q
\end{pmatrix}=
\begin{pmatrix}
1& uQ^{-1}\\ 0&1
\end{pmatrix}
\begin{pmatrix}
z-p&0\\
0&-Q
\end{pmatrix}
\begin{pmatrix}
1&0\\
Q^{-1}v& 1
\end{pmatrix}.
$$
One should now note that matrices of the form
\begin{tiny}
${\begin{pmatrix}
1&0\\
a&1
\end{pmatrix}}$
\end{tiny}
are always invertible with
\begin{tiny}
${\begin{pmatrix}
1&0\\
a&1
\end{pmatrix}^{-1}}=
{\begin{pmatrix}
1&0\\-a&1
\end{pmatrix}}$.
\end{tiny}
Thus the above calculation shows that $z-p$ is invertible if and only if $b-\hat p$ is invertible. Moreover, the inverses are related as follows:
$$
(b-\hat p)^{-1}
=\begin{pmatrix}
1&0\\
-Q^{-1}v& 1
\end{pmatrix}
\begin{pmatrix}
(z-p)^{-1}&0\\
0&-Q^{-1}
\end{pmatrix}
\begin{pmatrix}
1& -uQ^{-1}\\ 0&1
\end{pmatrix}
=\begin{pmatrix}
(z-p)^{-1}&\cdots\\
\cdots&\cdots
\end{pmatrix},
$$
and so we can get 
$G_p(z)=\ff((z-p)^{-1})$ as the (1,1)-entry of the matrix-valued
Cauchy-transform
$$G_{\hat p}(b)=\id\otimes \ff((b-\hat p)^{-1})=
\begin{pmatrix}
\ff((z-p)^{-1})&\cdots\\
\cdots&\cdots
\end{pmatrix}.
$$

We consider again the polynomial $p=xy+yx+x^2$ of our last example. Its selfadjoint linearization
can be written in the form
$$\hat p=
\begin{pmatrix}
0&0&0\\0&0&-1\\0&-1&0
\end{pmatrix}\otimes 1+
\begin{pmatrix}
0&1&\frac 12\\
1&0&0\\
\frac 12&0&0
\end{pmatrix} \otimes x+
\begin{pmatrix}
0&0&1\\0&0&0\\1&0&0
\end{pmatrix}\otimes y.$$
It is a linear polynomial in free variables, but with matrix-valued coefficients, and we 
need to calculate its matrix-valued Cauchy transform
$G_{\hat p}(b)=\id\otimes \ff((b-\hat p)^{-1})$.
This leads to the question if there exists a suitable matrix-valued version of free probability theory, with respect to 
the matrix-valued conditional expectation
$E=\id\otimes \ff$.

\section{Operator-valued extension of free probability}
\subsection{Basic definitions}
An operator-valued generalization of free probability theory was
provided by Voiculescu from the very beginning in \cite{V1985, V1995}.
The idea is that we replace our expectations, which yield numbers in $\CC$,
by conditional expectations, which take values in a fixed subalgebra $\cB$.
This is the analogue of taking conditional expectations with respect to sub-$\sigma$-algebras in classical probability.
Let us also remark that the concept of (operator-valued) freeness is distinguished on a conceptual level by symmetry considerations. In the same way as the classical de Finetti theorem equates conditionally independent and identically distributed random variables with random variables whose joint distribution
is invariant under permutations, a recent non-commutative version by K\"ostler and myself \cite{KSp} shows that in the non-commutative world 
one gets a corresponding statement by replacing ``conditionally independent'' by ``free with amalgamation'' and ``permutations'' by ``quantum permutations''. This has triggered quite some investigations on more general quantum symmetries and its relations to de Finetti theorems, see \cite{BCS, BaSp, RW, Web}.

\begin{definition}
1) Let $\cB\subset \cA$ be a unital subalgebra. A linear map
$E:\cA\to\cB$
is a \emph{conditional expectation} if
$E[b]=b$ for all  $b\in\cB$
and
$E[b_1ab_2]=b_1E[a]b_2$ for all $a\in\cA$ and $b_1,b_2\in\cB$.

2) An \emph{operator-valued probability space} consists of
$\cB\subset \cA$ and a conditional expectation $E:\cA\to\cB$.
If in addition $\cA$ is a
$C^*$-algebra, $\cB$ is a $C^*$-subalgebra of $\cA$, and $E$ is completely
positive, then we have an \emph{operator-valued $C^*$-probability space}.
\end{definition}

\begin{example}
Let $(\cA,\ff)$ be a non-commutative probability space. Put
$$M_2(\cA):=\left\{ \begin{pmatrix} a&b\\c&d\end{pmatrix}\mid a,b,c,d\in\cA\right\}$$
and consider $\psi:=\tr\otimes\ff$ and $E:=\id\otimes \ff$, i.e.:

$$\psi \left[\begin{pmatrix} a&b\\c&d\end{pmatrix}\right]=\frac {\ff(a)+\ff(d)}2,\qquad
E\left[ \begin{pmatrix} a&b\\c&d\end{pmatrix}\right]= \begin{pmatrix} \ff(a)&\ff(b)\\\ff(c)&\ff(d)\end{pmatrix}.$$
Then $(M_2(\cA),\psi)$ is a non-commutative probability space, and
$(M_2(\cA),E)$ is an $M_2(\CC)$-valued probability space.
\end{example}

Of course, we should also have a notion of distribution and freeness in the operator-valued sense.

\begin{definition}
Consider an operator-valued probability space $(\cA,E:\cA\to\cB)$. 
\begin{enumerate}
\item
The \emph{operator-valued distribution} of $a\in\cA$ is given
by all operator-valued moments
$E[ab_1ab_2\cdots b_{n-1}a]\in\cB~~(n\in\NN, b_1,\dots,b_{n-1}\in\cB)$.
\item
Random variables $x_i\in\cA$ ($i\in I)$ are \emph{ free with respect
to $E$} or \emph{free (with amalgamation) over $\cB$} if
$E[a_1\cdots a_n]=0$
whenever $a_i\in \cB\la x_{j(i)}\ra$ are polynomials in some $x_{j(i)}$ with coefficients from $\cB$,
$E[a_i]=0$ for all $i$, and $j(1)\not=j(2)\not=\dots\not=j(n).$
\end{enumerate}
\end{definition}

\begin{remark}
Polynomials in $x$ with coefficients from $\cB$ are linear combinations of expressions of the form $b_0xb_1xb_2\cdots b_{n-1}xb_n$.
It is important to note that the ``scalars'' $b\in\cB$ do not commute in general with the random variables $x\in\cA$.
\end{remark}

One can see that operator-valued freeness works mostly like ordinary freeness, one only has to take care
of the order of the variables. This means in all expressions they have to appear in their original
order.

\begin{example}
1) Note that in scalar-valued free probability one has the rule
$$\ff(x_1yx_2)=\ff(x_1x_2)\ff(y)\qquad\text{if $\{x_1,x_2\}$ and $y$ are free}.$$
By iteration, this
leads to a simple factorization of all ``non-crossing" moments in free variables. For example, if $x_1,\dots,x_5$ are free, then we have for the moment corresponding to
\setlength{\unitlength}{.4cm} $$\begin{picture}(9,4)\thicklines \put(0,0){\line(0,1){3}}
\put(0,0){\line(1,0){9}} \put(9,0){\line(0,1){3}} \put(1,1){\line(0,1){2}}
\put(1,1){\line(1,0){7}} \put(4,1){\line(0,1){2}} \put(8,1){\line(0,1){2}}
\put(2,2){\line(0,1){1}} \put(2,2){\line(1,0){1}} \put(3,2){\line(0,1){1}}
\put(5,2){\line(0,1){1}} \put(6,2){\line(0,1){1}} \put(6,2){\line(1,0){1}}
\put(7,2){\line(0,1){1}} \put(-0.1,3.3){$x_1$} \put(0.9,3.3){$x_2$} \put(1.9,3.3){$x_3$}
\put(2.9,3.3){$x_3$} \put(3.9,3.3){$x_2$} \put(4.9,3.3){$x_4$}\put(5.9,3.3){$x_5$} \put(6.9,3.3){$x_5$}
\put(7.9,3.3){$x_2$} \put(8.7,3.3){$x_1$}
\end{picture}$$
the factorization
$$
\ff(x_1x_2x_3x_3x_2x_4x_5x_5x_2x_1)=\ff(x_1 x_{1})\cdot
\ff(x_2 x_2 x_2)\cdot
\ff(x_3x_3)\cdot\ff(x_4)\cdot\ff(x_5x_5).
$$
This is actually the same as for independent classical random variables. The difference between classical and free shows up only for ``crossing moments".
\\
In the operator-valued setting one has the same factorizations of all non-crossing moments in free variables; but now one has to respect the order of the variables, the final expression
is of a nested form, corresponding to the nesting of the non-crossing partition. Here is
the operator-valued version of the above example.
\begin{eqnarray*}
E[x_1x_2x_3x_3x_2x_4x_5x_5x_2x_1]=E\Bigl[x_1\cdot E\bigl[x_2\cdot E[x_3x_3]\cdot x_2
\cdot E[x_4]\cdot E[x_5x_5]\cdot x_2\bigr]\cdot x_{1}\Bigr]
\end{eqnarray*}

2) For ``crossing" moments one also has analogous formulas as in the scalar-valued
case. But again one has to take care to respect the order of the variables. For example,
the formula
\begin{align*}
\ff(x_1x_2x_1x_2)=&\ff(x_1x_1)\ff(x_2)\ff(x_2)+\ff(x_1)\ff(x_1)\ff(x_2x_2)-\ff(x_1)\ff(x_2)\ff(x_1)\ff(x_2)
\end{align*}
for free $x_1$ and $x_2$ has now to be written as
\begin{multline*}
E[x_1x_2x_1x_2]=E\bigl[x_1E[x_2]x_1\bigr]\cdot E[x_2]+E[x_1]\cdot
E\bigl[x_2E[x_1]x_2\bigr]\\-E[x_1]\cdot E[x_2]\cdot E[x_1]\cdot E[x_2].
\end{multline*}
\end{example}
We see that, unlike in the scalar-valued theory, the freeness property in the operator-valued case uses the full nested structure of non-crossing partitions.

\subsection{Freeness and matrices}\label{sect:fm}
It is an easy but crucial fact that freeness is compatible with going over to matrices. For example
if 
$\{a_1,b_1,c_1,d_1\}$ and $\{a_2,b_2,c_2,d_2\}$ are free in $(\cA,\ff)$, then
\begin{tiny} ${\begin{pmatrix} a_1&b_1\\c_1&d_1\end{pmatrix}}$ \end{tiny} 
and \begin{tiny} ${\begin{pmatrix} a_2&b_2\\c_2&d_2\end{pmatrix}}$
\end{tiny} are in general not free in the scalar-valued probability space $(M_2(\cA),\tr\otimes\ff)$, 
but they are free with amalgamation over $M_2(\CC)$ in the operator-valued probability space $(M_2(\cA),\id\otimes\ff)$.

\begin{example}
Let $\{a_1,b_1,c_1,d_1\}$ and $\{a_2,b_2,c_2,d_2\}$ be free in $(\cA,\ff)$, consider
$$X_1:=\begin{pmatrix} a_1&b_1\\c_1&d_1\end{pmatrix}\qquad\text{and}\qquad
X_2:=\begin{pmatrix} a_2&b_2\\c_2&d_2\end{pmatrix}.$$
Then
$$X_1 X_2=\begin{pmatrix}
a_1a_2+b_1c_2&a_1b_2+b_1d_2\\c_1a_2+d_1c_2&c_1b_2+d_1d_2
\end{pmatrix}
$$
and for $\psi=\tr\otimes\ff$
\begin{align*}
\psi(X_1X_2)&=\bigl(\ff(a_1)\ff(a_2)+\ff(b_1)\ff(c_2)
+\ff(c_1)\ff(b_2)+\ff(d_1)\ff(d_2)\bigr)/2\\
&\not= (\ff(a_1)+\ff(d_1))(\ff(a_2)+\ff(d_2))/ 4\\
&=\psi(X_1)\cdot \psi(X_2),
\end{align*}
but for $E=\id\otimes\ff$
\begin{align*}
E(X_1X_2)&=\begin{pmatrix}
\ff(a_1a_2+b_1c_2)&\ff(a_1b_2+b_1d_2)\\
\ff(c_1a_2+d_1c_2)&\ff(c_1b_2+d_1d_2)
\end{pmatrix}
=
\begin{pmatrix} \ff(a_1)&\ff(b_1)\\ \ff(c_1)&\ff(d_1)\end{pmatrix}
\begin{pmatrix} \ff(a_2)&\ff(b_2)\\ \ff(c_2)&\ff(d_2)\end{pmatrix}\\
&=
E(X_1)\cdot E(X_2).
\end{align*}
\end{example}
Note that there is no comparable classical statement. Matrices of independent random variables do not show any reasonable structure, not even in an ``operator-valued" or ``conditional" sense.

\subsection{Operator-valued free cumulants}
In \cite{RS2} it was shown that the combinatorial description of free probability theory in terms of free cumulants can also be extended to the operator-valued setting. 
\begin{definition}
Consider $E:\cA\to\cB$. We define the \emph{free cumulants} $\kk_n^\cB:\cA^{n}\to\cB$ by
$$E[a_1\cdots a_n]=\sum_{\pi\in NC(n)} \kk_\pi^\cB [a_1,\dots,a_n].$$
The arguments of $\kk_\pi^\cB$ are distributed according to the blocks of $\pi$.
But now the cumulants are also nested inside each other according to the nesting of the blocks of $\pi$.
\end{definition}

\begin{example}
We consider 
$\pi=\bigl\{\{1,10\},\{2,5,9\},\{ 3,4 \} , \{ 6 \} , \{ 7,8 \} \bigr\}\in NC(10):$
\setlength{\unitlength}{.5cm} $$\begin{picture}(9,4)\thicklines \put(0,0){\line(0,1){3}}
\put(0,0){\line(1,0){9}} \put(9,0){\line(0,1){3}} \put(1,1){\line(0,1){2}}
\put(1,1){\line(1,0){7}} \put(4,1){\line(0,1){2}} \put(8,1){\line(0,1){2}}
\put(2,2){\line(0,1){1}} \put(2,2){\line(1,0){1}} \put(3,2){\line(0,1){1}}
\put(5,2){\line(0,1){1}} \put(6,2){\line(0,1){1}} \put(6,2){\line(1,0){1}}
\put(7,2){\line(0,1){1}} \put(-0.1,3.3){$a_1$} \put(0.9,3.3){$a_2$} \put(1.9,3.3){$a_3$}
\put(2.9,3.3){$a_4$} \put(3.9,3.3){$a_5$} \put(4.9,3.3){$a_6$} \put(5.9,3.3){$a_7$} \put(6.9,3.3){$a_8$}
\put(7.9,3.3){$a_9$} \put(8.7,3.3){$a_{10}$}
\end{picture}$$
For this we have
$$\kk^\cB_\pi[a_1,\dots,a_{10}]=\kk^\cB_2\Bigl(a_1\cdot
\kk^\cB_3\bigl(a_2\cdot\kk^\cB_2(a_3,a_4),a_5
\cdot\kk^\cB_1(a_6)\cdot\kk^\cB_2(a_7,a_8),a_9\bigr),a_{10}\Bigr).$$
\end{example}

\subsection{Operator-valued Cauchy and R-transform}
Now we consider operator-valued analogues of the Cauchy and $R$-transform. Again, those were introduced by Voiculescu, but without having the combinatorial meaning for the coefficients of the $R$-transform.
\begin{definition}
For $a\in\cA$, we define its \emph{operator-valued Cauchy
transform}
$$G_a(b):=E[\frac 1{b-a}]=\sum_{n\geq 0} E[b^{-1}(ab^{-1})^n]$$
and \emph{operator-valued $R$-transform}
\begin{align*}
R_a(b):&=\sum_{n\geq 0}
\kk_{n+1}^\cB(ab,ab,\dots,ab,a)=\kk_1^\cB(a)+\kk_2^\cB(ab,a)+\kk_3^\cB(ab,ab,a)+\cdots
\end{align*}
\end{definition}

As in the scalar-valued case we get as a relation between those two:
$$b G(b)=1+R(G(b))\cdot G(b)\qquad \text{or equivalently}\qquad G(b)=\frac 1{b-R(G(b))}.$$
If one reconsiders the combinatorial proof of these statements from the scalar-valued case, one notices that it respects the nesting of the blocks, so it works also in the operator-valued case. 

If one treats these concepts on the level of formal power series one gets all the main results as in the scalar-valued case, see \cite{Biane1, RS2, V1995, V2000}.

\begin{theorem}
If $x$ and $y$ are free over $\cB$, then:
mixed $\cB$-valued cumulants in $x$ and $y$ vanish;
it holds that $R_{x+y}(b)=R_x(b)+R_y(b)$;
we have the subordination
$G_{x+y}(z)=G_x(\omega(z))$.
\end{theorem}

\subsection{Free analysis}
In the last section we introduced the operator-valued $R$-transform and Cauchy transform on the level of formal power series.
In order to use them in an efficient way, we want to look at these objects in a more analytical way. This leads to the theory of ``free analysis".
This subject aims at developping a non-commutative generalization of holomorphic functions in the setting of operator-valued variables (or in the setting of several variables with the highest degree of non-commutativity). Free analysis was started by Voiculescu in the context of free probability
around 2000 \cite{V2000, FreeAnalysis1, FreeAnalysis2}; it builds on the seminal work of J.L.~Taylor \cite{Taylor}. Similar ideas are also used in work of Helton, Vinnikov and collaborators
around non-commutative convexity, linear matrix inequalities, or descriptor systems in
electrical engineering, see, e.g., \cite{KV}.

\subsection{Subordination in the operator-valued case}
Even more as in the scalar-valued theory it is hard to deal with the operator-valued $R$-transform in an analytical way.
Also, the operator-valued equation
$G(b)=\frac 1{b-R(G(b))}$ has hardly ever explicit solutions and, from the numerical point of view, it becomes quite intractable: instead of one algebraic equation we have now a system of algebraic equations.
However, there is also a
subordination version for the operator-valued case which was treated by Biane \cite{Biane1} and, more conceptually, by Voiculescu \cite{V2000}.

The following theorem shows that the analytic properties of the subordination function in the operator-valued situation are as nice as in the scalar-valued case.

\begin{theorem}[Belinschi, Mai, Speicher \cite{BMS}]\label{thm:main}
Let $(\cA, E:\cA\to\cB)$ be an operator-valued $C^*$-probability space and let
$x$ and $y$ be selfadjoint operator-valued random variables in $\cA$ which are
free over $\cB$. Then there exists a Fr\'echet analytic map $\omega\colon
\mathbb H^+(\cB)\to\mathbb H^+(\cB)$ so that
$$G_{x+y}(b)=G_x(\omega(b))\text{  for all } b\in\mathbb H^+(\cB).$$
Moreover, if $b\in\mathbb H^+(\cB)$, then $\omega(b)$ is the unique fixed point of the map
$$
f_b\colon\mathbb H^+(\cB)\to\mathbb H^+(\cB),\quad f_b(w)=h_y(h_x(w)+b)+b,
$$
and 
$$\omega(b)=\lim_{n\to\infty}f_b^{\circ n}(w) \qquad\text{for any $w\in\mathbb H^+(\cB)$}.$$ 
Here,
$\mathbb H^+(\cB):=\{b\in \cB\mid (b-b^*)/(2i)>0\}$ denotes the operator-valued upper halfplane of $\cB$, $h(b):=\frac 1{G(b)}-b$, and $b_b^{\circ n}$ is the $n$-th composition power of $f_b$ 
\end{theorem}

A similar description for the product of free variables in the operator-valued setting was shown by Belinschi, Speicher, Treilhard, and Vargas in \cite{BSTV}.

\section{Polynomials of independent random
matrices and polynomials in free variables}

Now we are able to solve the problem of calculating the distribution of a polynomial $p$ in free variables (and thus also the limiting eigenvalue distribution of the polynomial in asymptotically free random matrices).
The idea is to linearize the polynomial and to use operator-valued convolution for the linearization $\hat p$. We only present this for our running example.
The general case works in the same way.

\begin{example}
A linearization $\hat p$ of
$p=xy+yx+x^2$ was given in Eq.~\eqref{eq:phat} 
As we pointed out there, 
this means that the Cauchy transform $G_p(z)$ is given 
as the (1,1)-entry of the $M_3(\CC)$-valued Cauchy transform of $\hat p$:
$$G_{\hat p}(b)=\id\otimes\ff\left[(b-\hat p)^{-1}
\right]=
\begin{pmatrix}
G_p(z)&\cdots&\cdots\\
\cdots&\cdots&\cdots\\
\cdots&\cdots&\cdots
\end{pmatrix}
\qquad
\text{for} 
\quad
b=\begin{pmatrix}
z&0&0\\
0&0&0\\
0&0&0
\end{pmatrix}.$$
But now we can write $\hat p$ as $\hat p=\hat x+\hat y$ with selfadjoint
$$\hat x=
\begin{pmatrix}
0&\quad&x\quad&&\frac x2\\
x&\quad&0\quad&&0\\
\frac x2&\quad&0\quad&&0
\end{pmatrix}
\quad\text{and}\quad
\hat y=
\begin{pmatrix}
0&\quad&0\quad&&y\\
0&\quad&0\quad&&-1\\
y&\quad&-1\quad&&0
\end{pmatrix}
.$$
According to Sect.~\ref{sect:fm}, $\hat x$ and $\hat y$ are free over 
$M_3(\CC)$. Furthermore, the distribution of $x$ determines the operator-valued distribution of $\hat x$ and the distribution of $y$ determines the
operator-valued distribution of $\hat y$. This gives us the operator-valued Cauchy transforms of $\hat x$ and of $\hat y$ as inputs and we can use our
results on
operator-valued free convolution, in the form of Theorem \ref{thm:main}, to calculate the
operator-valued Cauchy transform of $\hat x+\hat y$ in
the subordination form
$G_{\hat p}(b)=G_{\hat x}(\omega(b))$,
where $\omega(b)$ is the unique fixed point in the upper half plane 
$\mathbb H^+(M_3(\CC))$ of the iteration
$$w\mapsto G_{\hat y}(b+G_{\hat x}(w)^{-1}-w)^{-1} - (G_{\hat x}(w)^{-1}-w).$$

There are no explicit solutions of those fixed point equations in $M_3(\CC)$, but a numerical implementation relying on iterations is straightforward.
One point to note is that $b$ as defined above is not in the open set $\mathbb H^+(M_3(\CC))$, but lies on its boundary. Thus, in order to be in the frame as needed in  Theorem \ref{thm:main}, one has to move inside the upper halfplane, by replacing 
$$b=\begin{pmatrix}
z&0&0\\
0&0&0\\
0&0&0
\end{pmatrix} \qquad \text{by}\qquad 
\begin{pmatrix} 
z&0&0\\
0&i\ee&0\\
0&0&i \ee
\end{pmatrix}$$
and send $\epsilon>0$ to zero at the end.

\begin{figure}[h]
\includegraphics[width=3.5in]{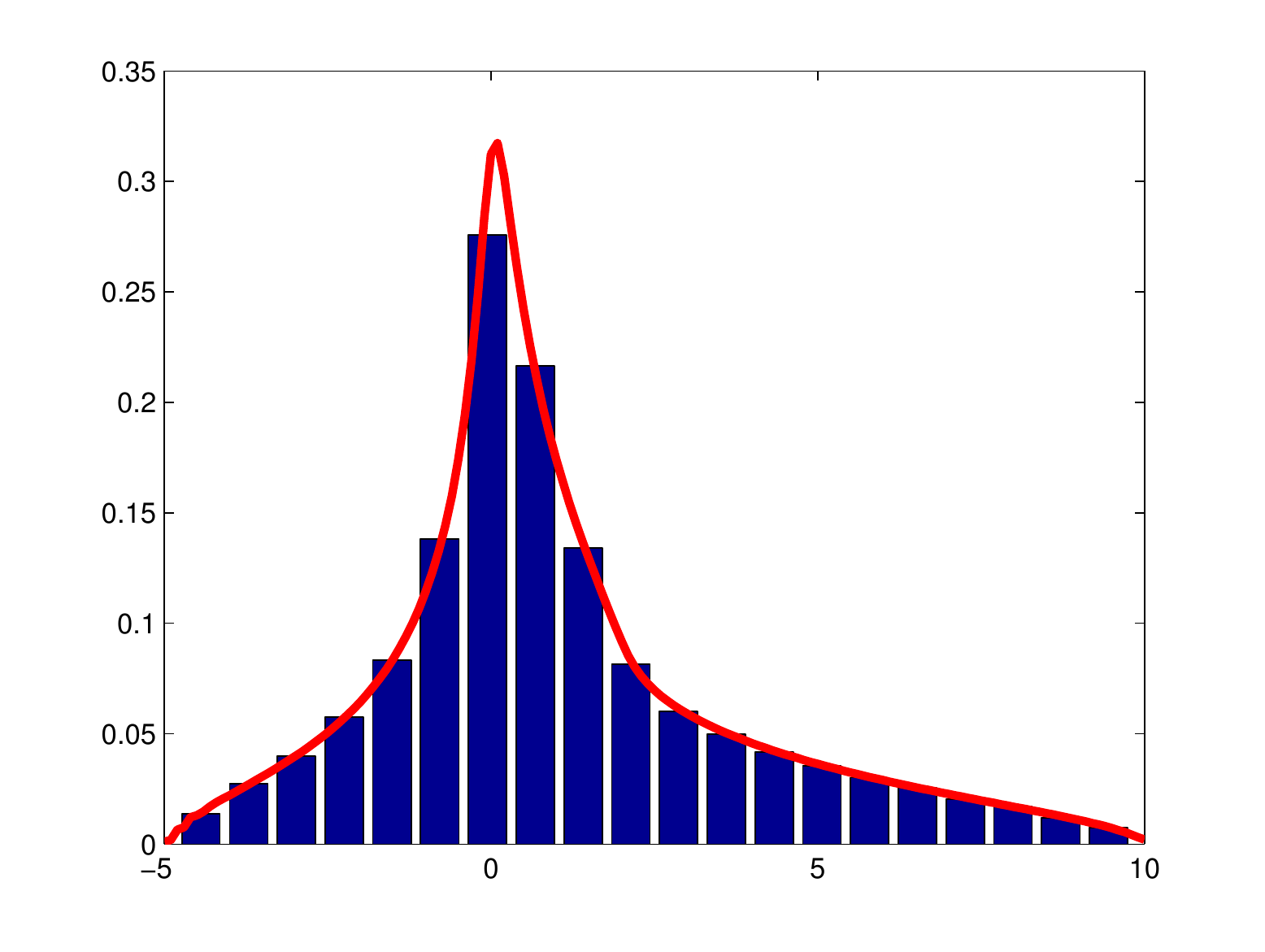}
\caption{\label{fig:final} Comparision between result of our algorithm for $p(x,y)=xy+yx+x^2$ ($x$ semicircular, $y$ Marchenko-Pastur) and histogram of eigenvalue distribution of $4000\times 4000$ random matrix $p(X,Y)$ (where $X$ and $Y$ are independent Gaussian and, respectively, Wishart random matrices)}
\end{figure}
\end{example}

\section{Further questions and outlook}\label{sect:outlook}
There are some canonical questions arising from this approach.

Firstly, our approach gives in principle a system of equations for 
the Cauchy transform of the wanted distribution. 
Whereas we can provide an efficient numerical fixed point algorithm for 
solving those equations, one would also like to derive qualitative properties of the solutions from this description. This will be pursued in the future. Prominent questions in this context are about the existence of atoms and regularity properties of the density of the distribution. One should note that,
by other approaches, Shlyakhtenko and Skoufranis made in \cite{ShSk} some
progress on such questions.

\begin{figure}[h]
$$\hskip-2cm
\includegraphics[width=3.5in]{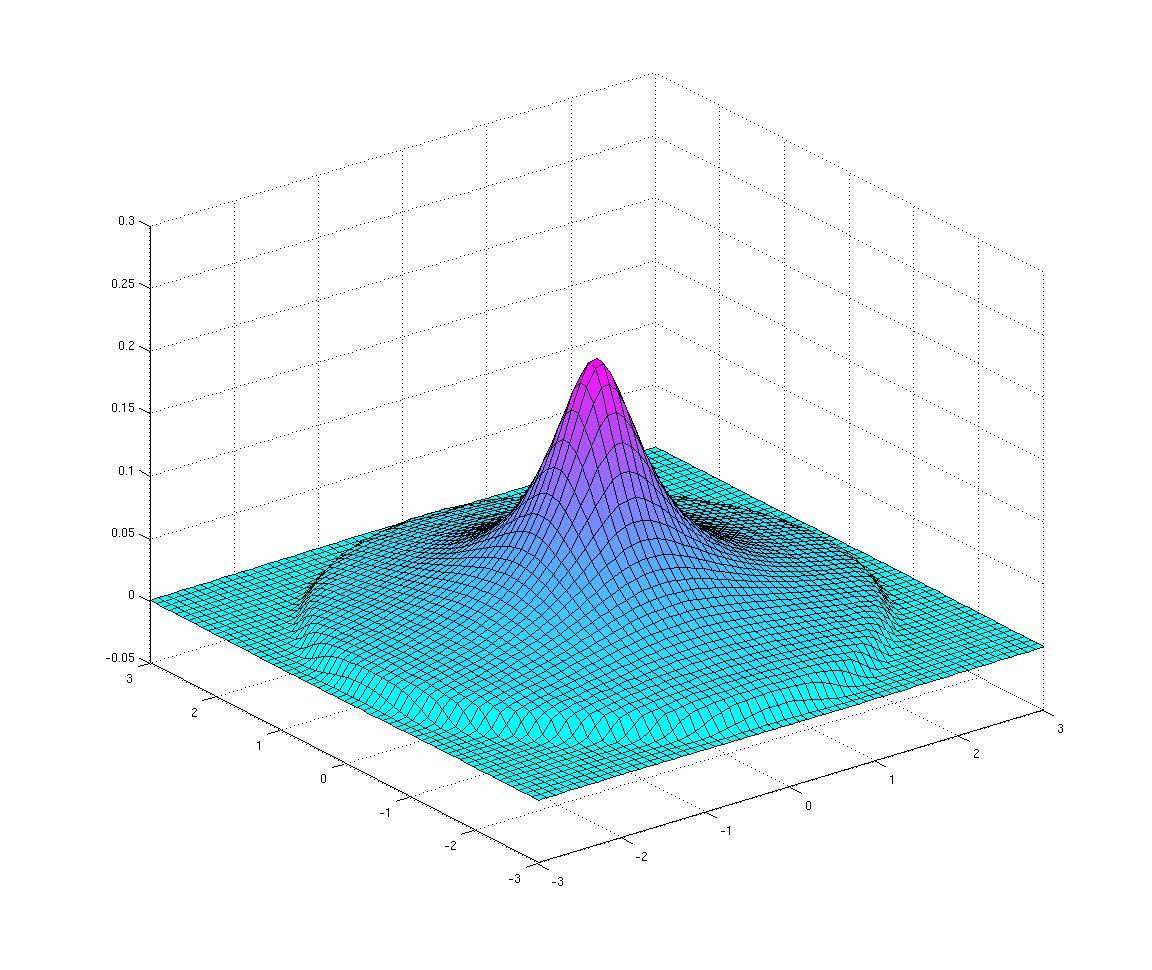}\hskip-1cm
\includegraphics[width=3.5in]{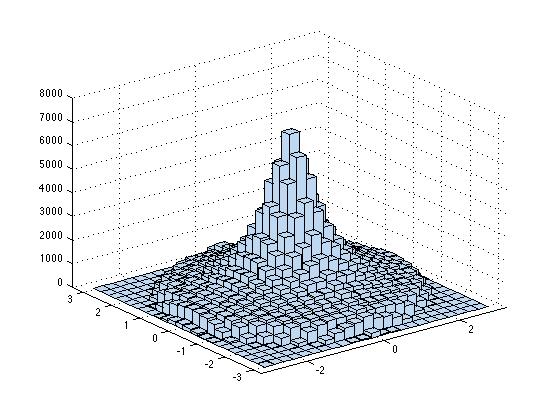}
$$
\caption{\label{fig:Brown} right: Brown measure, calculated by our algorithm, of the operator
$p(x,y,z) = xyz-2yzx+zxy$, where $x,y,z$ are free semicircular elements;
left: histogram of eigenvalues of $p(X_N,Y_n,Z_n)$, where $X_N,Y_N,Z_N$ are independent non-selfadjoint Gaussian $1000\times 1000$ random matrices; averaged over 100 realizations}
\end{figure}

Secondly, we concentrated here only on selfadjoint polynomials of selfadjoint variables, to ensure that we are dealing with selfadjoint operators. Then
the spectrum is a subset of the real line and thus the Cauchy transform contains all relevant information. In joint work with Belinschi and Sniday \cite{BSS} we are presently extending our ideas to non-selfadjoint polynomials, yielding non-normal operators. Then the spectral distribution of the operators has to be replaced by the so-called
Brown measure. By combining hermitian reduction ideas with the linearization trick and our subordination results one can then also extend our approach to this situation. 
An example for such a calculation is shown in Fig.~\ref{fig:Brown}.

Our methods should also work for more general classes of functions in non-commuting variables. In joint work with Mai, we are presently investigating the class of non-commutative rational functions.

\vskip-.0cm

\section*{Acknowledgement}
I am grateful to Alexandru Nica and the University of Waterloo
for their hospitality during the writing of this paper; and to Tobias Mai for preparing some of the figures. 
I thank Tobias Mai, Claus K\"ostler, and Moritz Weber for many useful comments on the first draft of this manuscript.

\end{document}